\documentclass[twoside,leqno]{article}

\usepackage[letterpaper]{geometry}

\usepackage{siamproceedings}

\usepackage[T1]{fontenc}
\usepackage{amsfonts}
\usepackage{graphicx}
\usepackage{epstopdf}
\usepackage{enumitem}
\usepackage{algorithmic}
\usepackage{amssymb}
\usepackage{amsmath}
\usepackage{array}
\usepackage{multirow}
\usepackage{mathrsfs}
\usepackage{mathtools}
\usepackage{esint}
\usepackage[normalem]{ulem}
\usepackage{marginnote}
\usepackage[english]{babel} 
\usepackage{threeparttable}
\usepackage{bm}
\usepackage{empheq}
\usepackage{algorithm}
\usepackage{algorithmic}

\usepackage{pgf,tikz}
\usetikzlibrary{arrows}
\usepackage[skip=2pt,font=scriptsize]{caption}
\usepackage{subcaption}

\usepackage{enumitem}

\usepackage{hyperref}
\usepackage{cleveref}
\usepackage{csquotes}
\hypersetup{
colorlinks,
linkcolor={red!80!black},
citecolor={green!80!black}
} 

\usepackage{makeidx}         
\usepackage{graphicx}        
\usepackage[multidot]{grffile}
\usepackage{multicol}        
\usepackage[bottom]{footmisc}

\ifpdf
  \DeclareGraphicsExtensions{.eps,.pdf,.png,.jpg}
\else
  \DeclareGraphicsExtensions{.eps}
\fi

\newsiamremark{remark}{Remark}
\newsiamremark{hypothesis}{Hypothesis}
\crefname{hypothesis}{Hypothesis}{Hypotheses}
\newsiamthm{claim}{Claim}

\DeclareMathOperator{\T}{\mathrm{T}}

\newcommand{\R}{\mathbb{R}}
\newcommand{\N}{\mathbb{N}}

\mathchardef\emptyset="001F

\renewcommand{\Sigma}{\boldsymbol{\sigma}}

\newcommand{\vecsigma}{\bm{\sigma}}
\newcommand{\zz}[2]{z_{#1}^{#2}}

\usepackage{amsopn}

\begin{document}

\newcommand\relatedversion{}

\title{\Large Machine Learning and Control: Foundations, Advances, and Perspectives
}
    \author{Enrique Zuazua\thanks{[1] Chair for Dynamics, Control, Machine Learning, and Numerics, Alexander von Humboldt-Professorship, Department of Mathematics,  Friedrich-Alexander-Universit\"at Erlangen-N\"urnberg, 91058 Erlangen, Germany; [2] Departamento de Matem\'{a}ticas, Universidad Aut\'{o}noma de Madrid, 28049 Madrid, Spain; [3] Chair of Computational Mathematics,  Deusto University, Av. de las Universidades, 24, 48007 Bilbao, Basque Country, Spain (\email{enrique.zuazua@fau.de}).}}

\date{}

\maketitle
 \fancyfoot[C]{\thepage}

\fancyfoot[R]{\scriptsize{Copyright \textcopyright\ 2026 by SIAM\\
Unauthorized reproduction of this article is prohibited}}

\begin{abstract} Control theory of dynamical systems offers a powerful framework for tackling challenges in deep neural networks and other machine learning architectures. \newline\indent
We show that concepts such as simultaneous and ensemble controllability offer new insights into the classification and representation properties of deep neural networks, while the control and optimization of static systems can be employed to better understand the performance of shallow networks. Inspired by the classical concept of {\it turnpike}, we also explore the relationship between dynamic and static neural networks, where depth is traded for width, and the role of transformers as mechanisms for accelerating classical neural network tasks. \newline\indent We also exploit the expressive power of 
 neural networks (exemplified, for instance, by the Universal Approximation Theorem) to develop a novel hybrid modeling methodology, the Hybrid-Cooperative Learning (HYCO), combining mechanics and data-driven methods in a game-theoretic setting. Finally, we describe how classical properties of diffusion processes, long established in the context of partial differential equations, contribute to explaining the success of modern generative artificial intelligence (AI).
\newline\indent We present an overview of our recent results in these areas, illustrating how control, machine learning, numerical analysis, and partial differential equations come together to motivate a fertile ground for future research. \end{abstract}

\section{Introduction.}\label{sec:1}
The interface between control theory (CT) and machine learning (ML), two disciplines with distinct foundations but increasingly convergent goals in intelligent systems, data-driven modeling, and scientific computing, is rapidly evolving. Control theory provides a rigorous foundation for feedback, stability, and optimization, principles that have long guided engineering systems. ML relies on data-driven optimization to model and predict complex, often unstructured phenomena.

In recent years, the boundaries between these two fields have become increasingly blurred \cite{zuazuanews}. Neural networks (NNs) can be viewed as discretized dynamical systems; training can be framed as an optimal control problem; and backpropagation mirrors sensitivity analysis in control. Likewise, when learning is data-driven but constrained by physical laws such as partial differential equations (PDEs), control insights become essential for designing stable and efficient models and algorithms.

These foundational links are not new. Aristotle already envisioned machines that could reduce human effort \cite{FCZ}. Centuries later, Wiener's \emph{Cybernetics} \cite{wiener1948cybernetics} gave this vision formal shape, defining it as ``the science of communication and control in animals and machines,'' and uniting control and communication.

We analyze classical ML models and NN architectures, including shallow NNs, deep residual networks (ResNets), neural ODEs (NODEs), and transformer architectures, together with fundamental questions such as data representation and generalization capacity. We show how control-theoretic ideas and methods provide new ways to address these challenges, yielding novel results, perspectives, and insights.

We conclude with related topics and future directions, focusing on the link between the classical theory of parabolic PDEs and generative diffusion models, federated learning, as well as the challenge of hybrid data-driven and physics-informed modeling.

By integrating control-theoretic thinking into ML, we deepen theoretical understanding while enhancing the reliability, interpretability, and efficiency of modern algorithms.

Further details can be found in the Oberwolfach Seminar Lecture Notes \cite{zuazuaober}.

\section{Shallow neural networks.}\label{SNN}

Shallow NNs constitute a mathematically tractable yet expressive framework for investigating foundational aspects of ML. Although structurally simple, they already encapsulate the essential challenges of representation, approximation, and generalization that persist in more sophisticated architectures. Their relative analytical accessibility renders them a natural model class for establishing rigorous connections between supervised learning practice and the theoretical methodologies of optimization, approximation theory, and numerical analysis.

\subsection{Supervised learning and data representation.} \label{SLUSL}

Supervised learning aims to learn a map \( f \colon \R^d \to \mathbb{R} \) from inputs (features) to outputs (labels) using labeled training data, constituting the observed finite dataset \( \{(x_i, y_i)\}_{i\in [N]} \subset \mathbb{R}^{d+1} \), where $[N]$ stands for the set of indices $[N]=\{ 1,..., N\}$. The goal is to minimize prediction error on unseen data, which inherently depends on how features capture the underlying structure of the problem. A good representation simplifies the input-output relationship, making it easier for the model to generalize. On the contrary, poor representations, in particular, those capturing noisy/irrelevant features, force the model to memorize unnecessary training data, leading to overfitting.

The origins of data representation trace back to the pioneering and visionary work of  Legendre (1805) and Gauss (1809) on the method of least squares, later extended using NN ans\"atze.

In line with \cite{liu2024representation}, we analyze shallow NN architectures within the framework of control and optimization.

\subsection{The model.}\label{modelNN}
We assume a consistent dataset, namely, a set of observations where all the input values (or features) $x_i \in \R^d$ are distinct. In practice, the observations may be corrupted by noise or measurement
errors. For the sake of simplicity, we consider scalar labels $y_i \in \R$, but our analysis can be easily extended to vector-valued labels.

The standard reconstruction methodology follows a least-squares approach, seeking the best approximation within a predefined function class specified by an NN ansatz, 
namely  a shallow NN involving $P$ neurons:
\begin{align}\label{eq:NN_shallow}
    f_{\text{shallow}} (x,\Theta) \coloneqq  \sum_{j\in [P]} w_j \sigma(\langle a_j, x \rangle + b_j), \, x \in \R^d.
\end{align}
Here $\sigma: \mathbb{R} \to \mathbb{R}$ denotes the activation function, 
typically the \emph{Rectified Linear Unit} (ReLU), i.e., $\sigma(s) = \max(0, s)= s_+$, and 
$\langle \cdot, \cdot \rangle$ is the standard inner product in $\mathbb{R}^d$. 
Furthermore, $P$ denotes the width of the NN, and each of the $P$ summands (neurons) depends on $d+2$ trainable parameters: the amplitude (or output weight) 
$w_j \in \mathbb{R}$, the bias (offset) $b_j \in \mathbb{R}$, and the input weights 
$a_j \in \mathbb{R}^d$. Collectively, these parameters are represented as
$
\Theta = \{ (w_j, a_j, b_j) \in \mathbb{R} \times \Omega \}_{j \in [P]},
$
where $\Omega$ is a compact subset of $\mathbb{R}^{d+1}$ containing a neighborhood of $0$. 

This model, commonly referred to in the literature as the two-layer NN, is the one that Cybenko employs in his pioneering paper on the universal approximation theorem (UAT) \cite{cybenko1989approximation}. His UAT result guarantees the density of this class of functions in the space of continuous functions on a hypercube of $\R^d$ for sigmoid activation functions $\sigma$, namely, a continuous function taking limits $0$ and $1$, as $s \to \pm \infty$.

Cybenko's result complements the pioneering work of Norbert Wiener in his renowned paper on Tauberian theorems \cite{wiener1932tauberian}. In \cite{wiener1932tauberian}, a simplified version of the ansatz \eqref{eq:NN_shallow} is considered: the activation function is a Gaussian $G$ defined in $\R^d$ (whose Fourier transform, a Gaussian as well, never vanishes),  the scaling parameter $a_j$ is omitted, and consequently, $b_j\in \R^{d}$. This leads  to a reduced formulation
 \begin{equation}\label{eq:NN_shallowb}
    f^G_{\text{shallow}} (x,\Theta) \coloneqq  \sum_{j\in [P]} w_j G(x  + b_j).
\end{equation}

In Cybenko's formulation (as in \eqref{eq:NN_shallow}), the Gaussian activation is replaced by a one-dimensional sigmoidal function designed to mimic the on-off behavior of biological neurons. To accommodate this change, Cybenko introduces a scaling factor $a_j$, which is not required in Wiener's one.

These two models are prototypes in nonlinear approximation theory, given that the ansatz depends not only on the amplitude parameters $w_j$, which enter linearly as multiplicative weights,  but also nonlinearly on the parameters $a_j$ and $b_j$. This is in contrast with the classical models of linear approximation theory, such as polynomial approximation (Stone--Weierstrass), Fourier series, wavelets, or the finite element method (FEM). 

Training these models, namely, determining the optimal parameter values $(w_j, a_j, b_j)$ (or $(w_j, b_j)$ in Wiener's model) that best represent the dataset by minimizing a loss function, leads to a nonconvex problem in the general setting of nonlinear approximation theory. The nonconvexity arises from the nonlinear dependence of the ansatz \eqref{eq:NN_shallow} on the parameters $(a_j, b_j)$. This section is primarily devoted to analyzing a convex relaxation of this minimization problem.


\subsection{Exact representation.}

The fact that a shallow NN can exactly interpolate the dataset whenever the number of neurons is greater than or equal to the cardinality of the dataset is a classical result, \cite{zhang2017understanding}. In particular, as shown in \cite{liu2024representation}, 
when $\sigma$ is continuous, $\sigma(s)=0$ for $s\leq 0$, and $\sigma(s)>0$ for $s>0$ (for instance the ReLU), and  $\Omega$ is a compact subset of $\R^{d+1}$ containing a neighborhood of $0$,
    for any consistent dataset $\{(x_i,y_i) \}_{i\in [N]} \subset \R^{d+1}$,
   and $P\geq N$,  there exists  $\Theta \in (\R\times \Omega)^{P}$ such that,  $f_{\textnormal{shallow}}$ as in \eqref{eq:NN_shallow} satisfies    \begin{equation}\label{exactrepresentation}
        f_{\textnormal{shallow}} (x_i,\Theta) = y_i, \, \forall i\in [N].
    \end{equation}

For a fixed dataset $\{(x_i,y_i) \}_{i\in [N]} \subset \R^{d+1}$, \eqref{exactrepresentation} guarantees the existence of parameters that yield an exact representation via~\eqref{eq:NN_shallow}, $P \ge  N$ being sufficient.  In some particular cases, exact representation might be achieved with fewer neurons. But, according to the theorem, $P = N$  suffices for all consistent datasets of $N$ pairs. 

The choice of the parameters $\Theta$ assuring exact representation is not unique. Indeed, generically with respect to $\{(a_j, b_j)\}_{j \in [P]}$, the basis functions involved are linearly independent, and the existence of weights $\{w_j\}_{j \in [P]}$ assuring \eqref{exactrepresentation} is then guaranteed. We are therefore in the regime known as overparameterized.

From a practical standpoint, it is not enough to know that a dataset can be represented; we also need effective methods to compute parameter values that ensure good generalization to unseen data.

\subsection{Optimal representation and relaxation.}\label{relaxation}

To address these issues, we adopt a complementary perspective by formulating an optimization problem that seeks parameter configurations which exactly represent the data while, among all such realizations, minimizing the $\ell_1$-norm of the neuron output weights $w_j$, namely:
\begin{equation}\tag{\text{P$_0$}}\label{pb:NN_exact}
  \begin{gathered}
    \inf_{\{(w_j,a_j,b_j)\in \R\times \Omega\}_{j\in [P]}}
      \sum_{j\in [P]} |w_j|,\\
    \text{s.t. }
      \sum_{j\in [P]} w_j\,\sigma\bigl(\langle a_j,x_i\rangle + b_j\bigr)
      = y_i, \, \forall i\in [N].
  \end{gathered}
\end{equation}
Note that only the norm of the weights $w_j$ is penalized in $\ell^1$, while the parameters $(a_j, b_j)$ are simply constrained to live in $\Omega$, a compact set of $\R^{d+1}$ containing the origin.

When the values \( \{y_i\}_{i\in [N]} \) represent noisy observations of the true underlying function, it becomes inappropriate to enforce exact interpolation. Instead, it is more meaningful to seek an approximate fit by relaxing the requirement of exact prediction. This results in an optimization problem formulated with an allowable margin of prediction error.

Specifically, for a prescribed error tolerance parameter \( \epsilon \geq 0 \), we consider the following optimization problem:
\begin{equation}\tag{\text{P$_\epsilon$}}\label{Pepsilon}
\begin{gathered}
  \inf_{\{(w_j,a_j,b_j)\in \R\times \Omega\}_{j\in [P]} } \sum_{j\in [P]} |w_j|,\\
  \text{s.t. } \left| \sum_{j\in [P]} w_j \sigma(\langle a_j, x_i \rangle + b_j) - y_i \right|
  \le \epsilon , \ \forall i\in [N].
\end{gathered}
\end{equation}
Problems \eqref{pb:NN_exact} and \eqref{Pepsilon} are nonconvex due to the nonlinearity of the activation function $\sigma$, and the nonlinear dependence on the parameters $(a_i, b_i)$, which induces the lack of convexity in their feasible sets. To cure this lack of convexity, we consider the following convex relaxation problems:
\begin{equation}\label{pb:NN_exact_rel}\tag{\text{PR$_0$}}
  \begin{gathered}
     \inf_{\mu\in \mathcal{M}(\Omega)} \|\mu\|_{\text{TV}}, \\
    \text{s.t. }  \int_{\Omega} \sigma( \langle a,x_i \rangle  + b)  d\mu(a,b) = y_i, \, \forall i\in [N];
  \end{gathered}
\end{equation}
and
\begin{equation}\label{pb:NN_epsilon_rel}\tag{\text{PR$_\epsilon$}}
  \begin{gathered}
     \inf_{\mu\in \mathcal{M}(\Omega)} \|\mu\|_{\text{TV}}, \\
     \text{s.t. }  \left|\int_{\Omega} \sigma( \langle a,x_i \rangle  + b)  d\mu(a,b) - y_i \right| \leq \epsilon , \, \forall i\in [N],
  \end{gathered}
\end{equation}
where  $\mathcal{M}(\Omega)$ represents the space of Radon measures on $\Omega$, and $\|\cdot\|_{\text{TV}}$ denotes the total variation norm.  These new relaxed optimization problems are formulated in the space of measures, under linear identity or inequality constraints, and they are clearly convex.

Before proceeding, we describe the link between the original NN ansatz $f_{\text{shallow}}$ and the one implicitly involved in this relaxed representation above, namely:
\begin{equation}\label{relaxedansatz}
f_{\text{relaxed}}(x, \mu) = \int_\Omega
\sigma( \langle a,x \rangle  + b)  d\mu(a,b). 
\end{equation}	
The function $f_{\text{shallow}}$ defined in \eqref{eq:NN_shallow} coincides with the representation in  \eqref{relaxedansatz} when the measure $\mu$ is atomic \begin{equation}\label{atoms}
\mu = \sum_{j\in [P]} w_j \delta_{(a_j, b_j)}.
\end{equation}
In fact, $f_{\text{relaxed}}$ in \eqref{relaxedansatz}  can be viewed as the continuous limit of \eqref{eq:NN_shallow}, obtained when the number of neurons tends to infinity and the discrete sum becomes a dense integral.

Given that exact representation is ensured with the primal discrete NN ansatz,  the relaxed ansatz \eqref{relaxedansatz} inherits this property. Accordingly, both the primal discrete problem (with $P \ge N$) and the convexified relaxation admit solutions. While the convex nature of the relaxed formulations suggests greater tractability, this advantage is offset by the fact that the relaxation operates in an infinite-dimensional space.

A natural question is whether the minimizers of the relaxed problem can be used to recover minimizers of the original finite-dimensional formulation. 

\subsection{On the lack of relaxation gap.}
In \cite{liu2024representation}, as a direct consequence of the representer theorem by Fisher-Jerome \cite{fisher1975spline},  we proved the following result guaranteeing that there is no gap between the primal problems and the relaxed ones, and that the extreme points of the relaxed solution sets have an atomic structure, with $N$ atoms, corresponding to the minimizers of the primal finite-dimensional optimization problem.

\begin{theorem}[No-Gap, \cite{liu2024representation}]\label{thm:NN_exists_0} 
  When $P \ge N$,
    the solution sets of \eqref{pb:NN_exact_rel} and \eqref{pb:NN_epsilon_rel}, denoted respectively by  $S\eqref{pb:NN_exact_rel}$ and  $S\eqref{pb:NN_epsilon_rel}$,  are nonempty, convex and compact in the weak-$*$ sense. 
    
    Moreover,
\begin{align}
  & \textnormal{val}(\ref{pb:NN_exact}) = \textnormal{val}(\ref{pb:NN_exact_rel}),
  && \textnormal{Ext}\bigl(S(\ref{pb:NN_exact_rel})\bigr)
     \subseteq
     \left\{ \sum_{j=1}^N w_j\delta_{(a_j,b_j)}
     \,\Big|\,
     (\Theta_j)_{j=1}^N \in S(\ref{pb:NN_exact})
     \right\}, \label{eq:exact}
  \\
  & \textnormal{val}(\ref{Pepsilon}) = \textnormal{val}(\ref{pb:NN_epsilon_rel}),
  && \textnormal{Ext}\bigl(S(\ref{pb:NN_epsilon_rel})\bigr)
     \subseteq
     \left\{ \sum_{j=1}^N w_j\delta_{(a_j,b_j)}
     \,\Big|\,
     (\Theta_j)_{j=1}^N \in S(\ref{Pepsilon})
     \right\}, \label{eq:epsilon}
\end{align}
    where $\textnormal{Ext}(S)$ represents the set of all extreme points of $S$, $S(P)$ the solution set of the corresponding optimization problem $P$ and $ \textnormal{val}(P)$ the minimum value.
\end{theorem}

In light of this result, the values of the primal problems become independent of $P$ as soon as $P \geq N$. Consequently, increasing $P$ beyond $N$ does not improve the optimal value of the objective function. At the same time, choosing $P = N$ minimizes the representational complexity of the NN ansatz, making it a natural and efficient hyperparameter choice. On the other hand, the result guarantees that the extremal minimizers of the relaxed problem have an atomic structure, leading to the original NN ansatz.

This result paves the way for efficient optimization methods by exploiting the convexity of the relaxed problem.

\subsection{Generalization error.}

In practical applications, the primary goal of supervised learning is to approximate an unknown target function. Consequently, NNs must demonstrate strong generalization capabilities, not only within the training dataset but also beyond it,  in the more challenging setting of \emph{out-of-distribution generalization}.

 We assess this via a testing dataset $\{(X', Y')\} = \{(x_i', y_i') \in \mathbb{R}^{d+1}\}_{i=1}^{N'}$, where $N' \in \mathbb{N}_+$, distinct from the training data, drawn independently from the same underlying distribution. 

The network's generalization quality is quantified by its performance on $(X', Y')$, measured through a comparison between the true outputs $\{y_i'\}_{i=1}^{N'}$, and
the corresponding predictions $\{f_{\mathrm{shallow}}(x_i', \Theta)\}_{i=1}^{N'}$.

There are several ways of proceeding to this comparison. Rather than focusing on pointwise errors, we analyze discrepancies at the distributional level to gain more robust insights. This approach provides several advantages: a more comprehensive evaluation of model behavior, reduced sensitivity to individual outliers, and better alignment with statistical learning theory.

Let us denote by 
\begin{equation}
     m_X = \frac{1}{N}\sum_{i\in [N]} \delta_{x_i}, \quad m_Y = \frac{1}{N}\sum_{i\in [N]} \delta_{y_i}, \quad      m_{\textnormal{train}}= \frac{1}{N}\sum_{i\in [N]} \delta_{(x_i,y_i)}, 
    \end{equation} 
    the empirical distribution of the training dataset, and by
    \begin{equation}
m_X' = \frac{1}{N'}\sum_{i=1}^{N'} \delta_{x_i'}, \quad m_Y' = \frac{1}{N'}\sum_{i=1}^{N'} \delta_{y_i'}, \quad    m_{\textnormal{test}}= \frac{1}{N'}\sum_{i=1}^{N'} \delta_{(x_i', y_i')}, 
    \end{equation} 
    the one of the testing datasets.
    Let us also consider the measure of predictions:
     \begin{equation}
          m_{\textnormal{pred}} (\Theta)= \frac{1}{N'}\sum_{i\in [N']} \delta_{(x_i', f_{\text{shallow}} (x_i', \Theta))}.
     \end{equation}
  
  In this measure theoretical setting, distances between the different datasets will be calibrated in terms of the  Kantorovich--Rubinstein (KR) distance $d_{\textnormal{KR}}(\cdot, \cdot) $, that respects the underlying geometry of the data, accounting for how far points must be transported to transform one dataset into another.
    It is equivalent to the $1$-Wasserstein distance, defined as follows:
    \[
    \mathcal{W}_1(\mu, \nu) \coloneqq  \min_{\pi \in \Pi(\mu, \nu)} \int_{\R^d \times \R^d} |x - y| \, d\pi(x, y),
    \]
    where \( \Pi(\mu, \nu) \) denotes the set of measures \( \pi \) on \( \R^d \times \R^d \) that couple \( \mu, \nu \in \mathcal{P}_{ac}^c(\R^d) \), two compactly supported probability measures on $\R^d$.

     \begin{theorem}[Generalization bound, \cite{liu2024representation}]\label{thm:generalization}   
    Assume that the activation function \(\sigma\) is \(\mathcal{L}\)-Lipschitz. Then, for any $\Theta\in \R^{(d+2)P}$, we have
    \begin{equation}\label{eq:generalization}
        d_{\textnormal{KR}}(m_{\textnormal{test}}, m_{\textnormal{pred}} (\Theta)) \leq \underbrace{d_{\textnormal{KR}}(m_{\textnormal{train}}, m_{\textnormal{test}})}_{\textnormal{Irreducible error from datasets}} + r(\Theta),
    \end{equation}
    where
    \begin{equation}\label{eq:r}
        r(\Theta) = \underbrace{\frac{1}{N} \sum_{i\in [N]} \left|f_{\textnormal{shallow}} (x_i, \Theta) - y_i\right|}_{\textnormal{Training bias}} + \underbrace{d_{\textnormal{KR}}(m_X, m_{X'}) \mathcal{L} \sum_{j\in [P]} |w_j|\|a_j\|}_{\textnormal{Sensitivity term}}.
    \end{equation}
    
    Furthermore, let \( P \geq N \), \(\sigma\) be the ReLU activation function (so that $\mathcal{L}=1$) and \(\Omega\) be the unit ball in \(\mathbb{R}^{d+1}\). Then, for any \(\epsilon \geq 0\), the solution \(\Theta_{\epsilon}\)  of \eqref{Pepsilon} satisfies
    \begin{equation}
        r(\Theta_{\epsilon}) \leq \mathcal{U}(\epsilon)\coloneqq \epsilon + d_{\textnormal{KR}} (m_X,m_{X'})\,\textnormal{val}\eqref{pb:NN_epsilon_rel}. \label{eq:generalization_epsilon}    \end{equation}
\end{theorem}

The irreducible error in \eqref{eq:generalization} is independent of the parameters \(\Theta\), being completely determined by the training and testing datasets. 
In contrast, the residual term \(r(\Theta)\) consists of two components depending on \( \Theta \): (1) the fidelity error $\epsilon$ on the training set, referred to as the training bias, and (2),  the sensitivity or stability of the trained NN.

The second part of the theorem provides a sharp bound on the term $r(\Theta) $ for the optimal choice of the parameters $\Theta_\epsilon$.
In light of the estimate \eqref{eq:generalization_epsilon}, the problem of minimizing the right-hand-side upper bound with respect to $\epsilon$ arises. The analysis of this problem is conducted in \cite{liu2024representation}.

The main conclusion is that two scenarios have to be distinguished:
\begin{itemize}
\item  The exact representation problem \eqref{pb:NN_exact} with \(\epsilon = 0\) is sufficient to guarantee good generalization properties, i.e., the NN can represent the test distribution without needing to introduce any approximation error \(\epsilon > 0\), when the distance
\( d_{\textnormal{KR}}(m_{\textnormal{train}}, m_{\textnormal{test}}) \) between the training and testing distributions lies below a certain threshold  
(which can be identified through dual problem analysis; see \cite{liu2024representation}). Moreover,  the generalization performance inevitably deteriorates as \(\epsilon\) and  \(\mathcal{U}(\epsilon)\) grow.

\item  Conversely, when 
\( d_{\textnormal{KR}}(m_{\textnormal{train}}, m_{\textnormal{test}}) \) exceeds this threshold, the optimal 
$\epsilon$
 becomes strictly positive and can be determined by solving the dual problem of \eqref{pb:NN_epsilon_rel} (see \cite{liu2024representation}). In this regime, choosing a smaller value of 
$\epsilon$ would harm generalization performance.
\end{itemize}

\subsection{Conclusions.}

This section summarizes the main insights of \cite{liu2024representation} on shallow NNs through the lens of static optimal control and convex optimization. Training with an $\ell_1$-penalty on output weights admits a convex relaxation in the space of measures with no relaxation gap: extreme points of the relaxed solution set are $N$-atomic, recovering the finite NN. Two hyperparameters follow from this structure:
\begin{itemize}
\item Width: Exact interpolation requires at most $P=N$; taking $P>N$ does not improve the optimal value.
\item Tolerance: The generalization-optimal tolerance $\epsilon$ in the relaxed problem depends on the distribution shift $d_{\textnormal{KR}}(m_{\textnormal{train}}, m_{\textnormal{test}}) $: below a data-dependent threshold, $\epsilon^\star=0$ (exact fit) is optimal; above it, $\epsilon^\star>0$ is determined by the dual of the relaxed program.
\end{itemize}
Thus, the convex formulation simultaneously identifies minimal representational width and prescribes the fidelity level needed for out-of-distribution robustness as a function of \( d_{\textnormal{KR}}(m_{\textnormal{train}}, m_{\textnormal{test}}) \).

For details, see \cite{liu2024representation}, which also develops efficient algorithms for the high-dimensional (infinite-dimensional in the relaxed setting) optimization problems, addressing the classical curse of dimensionality.

\section{Control-based supervised learning with residual neural networks and neural ODEs.}\label{sec:2}
Previously, we showed that shallow NNs can represent finite datasets, and we also derived estimates on the model complexity and on their generalization capacity. We further observed that relaxation allows us to convexify the training problem, at the price of making it infinite-dimensional.

Deeper insights arise from analyzing the gradient descent dynamics when training the finite-dimensional ansatz. This is a delicate matter, as nonlinear phenomena such as condensation may occur \cite{zhou2022towards}. These effects are linked to the absence of a gap between the original finite-dimensional training problem and its relaxed infinite-dimensional counterpart. Indeed, even when starting from a continuum of parameters, $\ell^1$-optimal minimizers involve only $N$ neurons. Condensation also accounts for the fact that exact representation can sometimes be achieved even when $P < N$, or even in the underparameterized regime $P < N/(d+1)$. However, a complete characterization of the datasets that can be exactly represented with $P < N$ neurons remains out of reach.

The methodology based on the static NN ansatz \eqref{eq:NN_shallow} offers only limited interpretative power. We therefore turn to deep NNs, where the role of width in shallow architectures is replaced by depth, modeled as a discrete iterative dynamical system. This viewpoint reformulates representation as a control problem, providing deeper insight into how parameters should be selected to perform data representation tasks effectively. It also leads to constructive strategies that not only identify suitable networks but also yield quantitative estimates. In turn, this framework offers valuable guidance on what can be expected when training deep networks through optimization methods, as is commonly done in practice.

\subsection{Supervised learning as a control problem.}

We focus on residual NNs (ResNets) introduced in \cite{he2016deep} (see also \cite{weinan2017proposal}, \cite{haber2017stable}, and \cite{li2018maximum}). A ResNet of depth \( L \in \mathbb{N} \) can be described as a discrete-time dynamical system evolving in $\R^d$:
\begin{align}\label{eq:ResNet}
    x^{k+1} = x^{k} + W^k \vecsigma(A^k x^k + b^k), \quad k \in [L],
\end{align}
where:
\begin{itemize}
    \item The index \( k \) acts as a pseudo-time variable in this discrete-time system.
    \item \( x^k \in \R^d \) represents the evolving state of the network.
    \item The parameters \( W^k \in \R^{d \times P} \), \( A^k \in \R^{P \times d} \), and \( b^k \in \R^P \) define the transformations at each layer.
    \item \( L \geq 1 \) and \( P \geq 1 \) denote the depth and width of the network, respectively.
    \item \( \vecsigma: \R^P \to \R^P \) is a vector-valued activation function, in which the scalar one $\sigma$ is applied component-wise. A common choice is the (vectorial) ReLU function.
\end{itemize}

When \( P = 1 \), the model simplifies to:
\begin{align}\label{eq:ResNetscalar}
    x^{k+1} = x^{k} + w^k \sigma(\langle a^k, x^k \rangle + b^k), \quad k \in [L],
\end{align}
where \( b^k \in \R \) are scalar drift coefficients and \( a^k, w^k \in \R^d \) are vector valued, respectively.

The ResNet architecture differs fundamentally from the shallow model \eqref{eq:NN_shallow}. Rather than stacking all neurons simultaneously, ResNets introduce them progressively (one in \eqref{eq:ResNetscalar} or \( P \) in \eqref{eq:ResNet}, at a time) through an inductive, layer-wise construction. This shift from width (represented by \( P \) in shallow networks) to a combination of depth  \( L \) and width $P$  in ResNets provides several advantages. Notably, as we shall see, the representational properties of ResNets can be naturally interpreted within a control-theoretical framework.

The ResNet \eqref{eq:ResNet} can be viewed as the Euler discretization of the continuous-time dynamics  known as  NODE:
\begin{equation}\label{eq:node}
   \dot{x}(t) = W(t) \vecsigma(A(t) x(t) + b(t)), \quad t \in (0, T); \quad 
        x(0) = x^0.
\end{equation}

This provides a continuous interpretation of deep learning dynamics, bridging discrete architectures with differential equations.
In this continuous-time model:
\begin{itemize}
    \item The state \( x(t) \in \R^d \) evolves continuously over time.
    \item The parameters \( (W(t), A(t), b(t)) \), which play the role of controls, depend on the continuous time $t$. We can assume them to lie in the space $\mathcal{S}_P=L^\infty\left(0, T; \R^{d \times P}\times \R^{P \times d} \times \R^P)\right).$ The controls we construct are piecewise constant with finitely many jumps, and therefore belong to the space of functions of bounded variation, denoted by $BV\left(0, T; \R^{d \times P}\times \R^{P \times d} \times \R^P\right)$.

    \item The time horizon $T > 0$ can be interpreted as the depth $L$ of the corresponding ResNet. A more accurate perspective, however, is to regard the depth $L$ as the number of control switches in the associated NODE when the controls are chosen to be piecewise constant, as we shall explain in detail below.
   
\end{itemize}

The NODE model defines the flow map:
\[
    \Phi_T(\cdot; W, A, b) : x^0 \in \R^d \mapsto x(T; x^0) \in \R^d,
\]
where \( x(t; x^0) \) is the solution to \eqref{eq:node} with initial condition \( x^0 \) and the chosen time-dependent controls.

The representation problem in supervised learning can now be reformulated in control-theoretic terms: the task becomes that of mapping the input data (considered as the initial data of the ResNet) to their corresponding labels (the targets at the final time $T$) via the flow map $\Phi_T(\cdot; W, A, b)$, through a suitable choice of the control parameters $(W(t), A(t), b(t))$ of the NODE.

A fundamental difference between this control-theoretic interpretation of supervised learning and classical control problems lies in the fact that, in supervised learning, the same controls $(W(t), A(t), b(t))$ must simultaneously drive the entire ensemble of input data to their corresponding outputs. This makes the problem one of ensemble, or simultaneous, control, unlike the classical setting, where controls are typically tailored to each specific initial state and target.

For a NODE to solve a representation problem, the dataset must satisfy a consistency condition: distinct initial states cannot map to the same output due to the forward and backward uniqueness property of ODEs and, in particular, of NODEs. In one spatial dimension (\(d=1\)), the flow map is monotonically increasing (for Lipschitz vector fields) because ODE trajectories cannot cross: larger initial conditions always yield larger states at any later time. This imposes a severe limitation on controllability, and in particular, the simultaneous control property discussed above fails. For this reason, we will focus on the case \(d \geq 2\), where NODEs become significantly more expressive. In higher dimensions, the monotonicity constraint of the one-dimensional case no longer applies, since trajectories may bend and overlap in projection. Thus, no order-preserving structure restricts the flow map. 
  
Here, inputs and outputs (or labels) both lie in the same space $\R^d$. In practical applications, however, this is generally not the case. For instance, in the shallow NN framework considered in the previous section, the inputs belonged to $\R^d$, while the outputs were in $\R$. NODEs can be naturally adapted to such situations by composing the flow map $\Phi_T(\cdot; W, A, b)$ with a suitable linear projection or extension operator that maps the ambient space $\R^d$ into the Euclidean space of labels. This is, in fact, the standard practice in most applications, also when dealing with NODEs. This issue will be omitted here for the sake of brevity.

\subsection{Simultaneous control: depth versus width.}\label{subsec:simultaneous_control}

A fundamental challenge in control theory is determining whether a system can be fully manipulated to achieve a desired outcome. This is formalized through the concept of controllability, which explores whether it is possible to steer a system from any given initial state to any desired final state within a finite time, using admissible controls.
However, as mentioned above, with supervised learning applications in mind, in the context of NODEs, we are rather interested in a more demanding goal: that of simultaneous or ensemble control.

We focus on piecewise-constant controls in time that, as we shall see, will suffice to achieve our goals. This choice, on the one hand, naturally induces a layered structure that mirrors the discrete nature of ResNet architectures. On the other hand, it renders the NODE dynamics more tractable and interpretable, while also facilitating the design of control inputs.  We will be able to prove the needed property of simultaneous control by carefully and inductively defining each of the values that the controls take. 

To break this down, we consider the equivalent formulation of \eqref{eq:node}:
\begin{equation}\label{eq:node-p}
    \dot{x}(t) = \sum_{i\in[P]} w_i(t)\sigma(\langle a_i(t),  x(t) \rangle + b_i(t)), \hspace{1cm} t \in (0,T),
\end{equation}
where $w_i(t), a_i(t)\in\R^d$ are the columns of $W(t)$ and the rows of $A(t)$, respectively, and $b_i(t) \in\R$ is the $i$-th coordinate of $b(t)$. 

In what follows, to streamline the presentation, we focus on the case where the activation function $\sigma$ is the ReLU, although most of what we say can be easily generalized for a broad class of activation functions.

For insight,  note that for each \(t \in (0,T)\) and \(i \in [P]\), the controls \(a_i(t)  \in \R^d\) and \(b_i(t)  \in \R\) define a $(d-1)$-dimensional hyperplane  
\[
H_i(t) := \{ x \in \R^d : \langle a_i(t), x \rangle + b_i(t)= 0 \},
\]  
that partitions the Euclidean space \(\R^d\) into two disjoint complementary half-spaces
\begin{equation}\label{splitting}
   \begin{array}{l}
   \begin{cases}
    H^+_i(t)\coloneqq\big\{x\in\R^d:\langle a_i(t),  x\rangle+b_i(t)>0\big\},\\ H^-_i(t)\coloneqq\R^d\setminus H^+_i(t)\coloneqq\big\{x\in\R^d:\langle a_i(t),  x\rangle+b_i(t)\le 0\big\}.
      \end{cases}  \end{array}
    \end{equation}
\begin{figure}[h]
    \centering
    \includegraphics[width=0.25\linewidth]{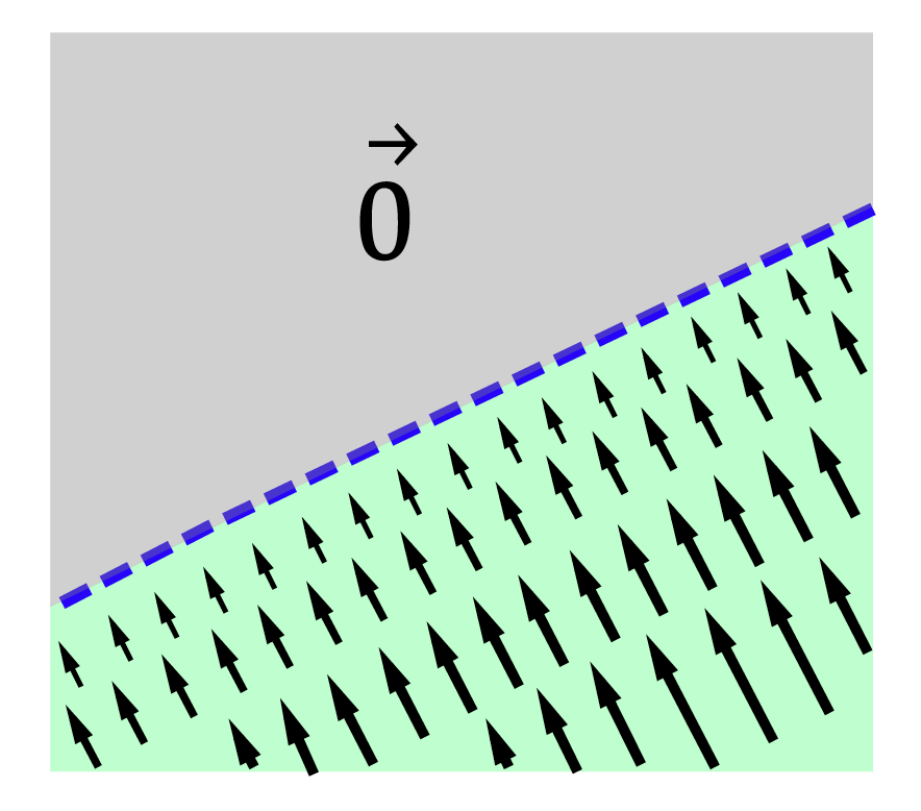}
    \includegraphics[width=0.25\linewidth]{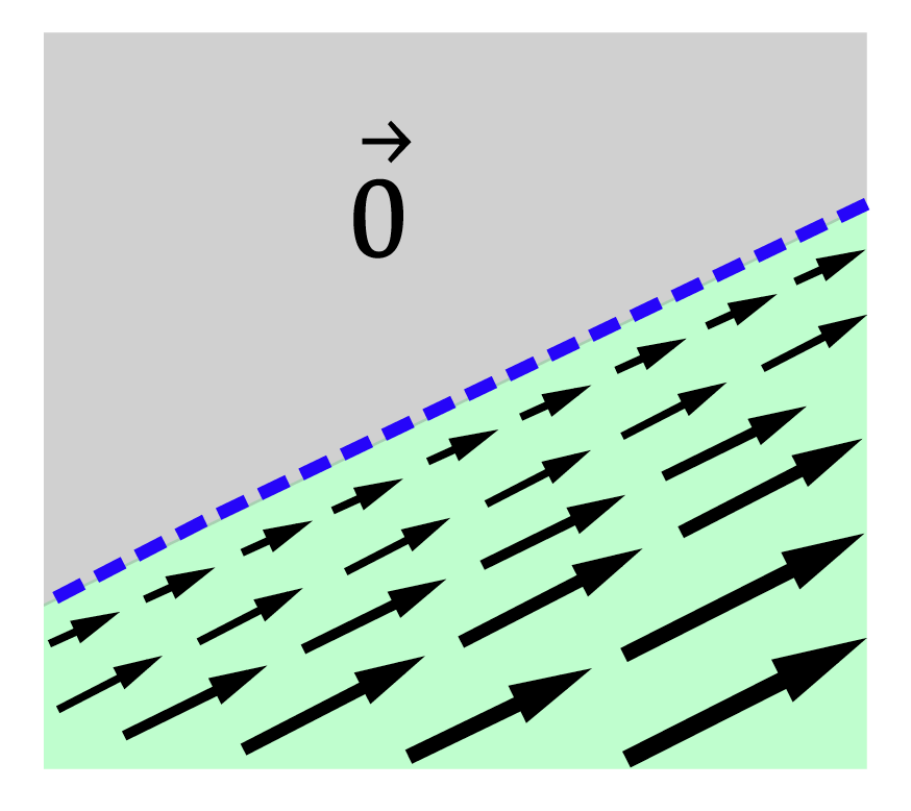}
    \includegraphics[width=0.25\linewidth]{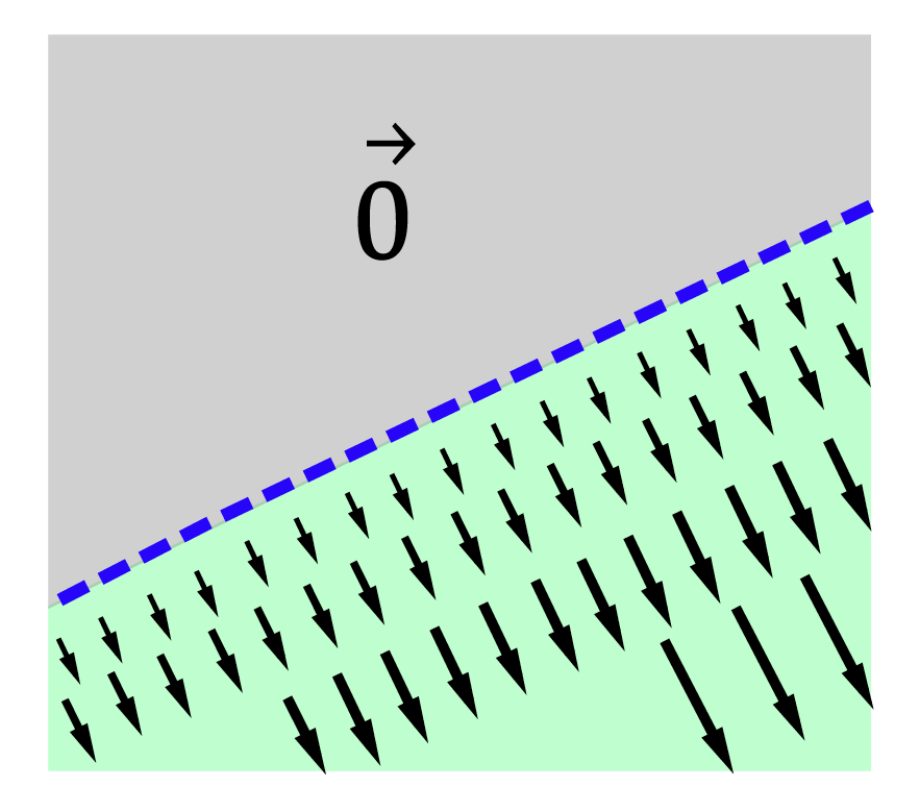}
    \caption{Basic movements generated by each of the neurons in the vector field determining the dynamics of \eqref{eq:node}: $(a_i, b_i)$ determines the cutting hyperplanes defining the various regions, and $w_i$ specifies the direction of the ``wind'' on the active half-space, namely, the side where the neuron drives the dynamics.
  The blue dashed line represents the hyperplane $ \langle a_i,  x\rangle+b_i=0$. From left to right: Compression, parallel motion, and expansion, depending on the choice of $w_i$. }
    \label{fig1}
\end{figure}

\begin{figure}[h]
    \centering
  \includegraphics[width=0.25\linewidth]{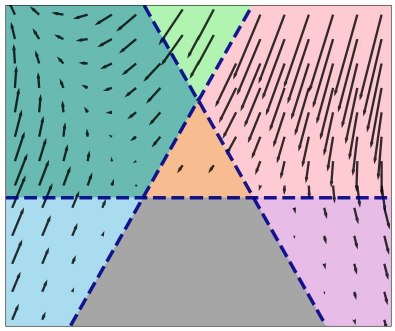} \hspace{1cm}
    \includegraphics[width=0.25\linewidth]{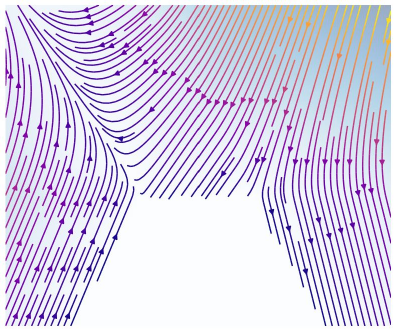}
 
      \caption{The Voronoi-type decomposition determined by the neural vector field at any time when involving several neurons $P \ge 2$. In the present case, $P=3$ neurons decompose the plane ($d=2$) into 7 regions where the neural vector field is oriented differently. In the figure on the right, we observe that the vector field vanishes within the lower-center Voronoi cell, thereby defining a region of space that remains stationary under the action of the vector field. }
    \label{Voronoi}
\end{figure}

Meanwhile, each control \( w_i(t)  \in \R^d\) acts only on the points within the half-space \(H^+_i(t)\), given that  $\sigma(\langle a_i(t),  x\rangle + b_i(t))=0$ for all $x\in H^-_i(t)$. This generates a linear motion on \(H^+_i(t)\) when $\sigma$ is the ReLU. Without loss of generality, by scaling,   the controls $a_i(t)$ are  normalized to unit norm, i.e., $|a_i(t)|=1$. By doing so, each neuron on the vector field determining the dynamics of the NODE takes the form $w_i(t)\operatorname{dist}( x, H^-_i(t))  $ (see Figure \ref{fig1}). Summing up the contribution of the $P$ neurons, the vector field driving the dynamics in \eqref{eq:node} is given by a weighted superposition of the form \(\sum  w_i(t) \operatorname{dist}( x, H^-_i(t)) \) on each \( x \in \R^d\), the sum being taken on the active neurons on $x$. This defines a piecewise linear and continuous vector field on the Voronoi-type decomposition induced by the hyperplanes $H_i(t)$, which partition the space into polytopes (see Figure \ref{Voronoi}).

The case $P=1$ was addressed in \cite{domenecNODES}, showing the aimed simultaneous controllability for consistent datasets. In this case the model reads: \begin{equation}\label{eq:nodep1}
    \dot{x}(t) = w(t) \, \sigma(\langle a(t),  x(t)\rangle + b(t)), \hspace{1cm} t \in (0, T).
\end{equation}
Here, at any given time  $t$,  only a single neuron is active. The vector field driving the dynamics vanishes on one half-space and is linear on the other. This structure allows us to design an inductive strategy that builds piecewise constant controls 
\((w, a, b)\), sequentially steering each input point \(x_i\) to its prescribed target 
\(y_i\),  thereby yielding the following result.

\begin{theorem}[Simultaneous controllability, $P=1$, \cite{domenecNODES}]\label{thm:1}
    Let the dimension $d \ge 2$, $N\geq1$, and $T>0$ be fixed. Given a consistent dataset \(\{ (x_i, y_i)\}_{i\in[N]}\subset\mathbb{R}^{d}\times \mathbb{R}^{d}\) with $x_i\neq x_j$ and $y_i\neq y_j$ if $i\neq j$, there exists  piecewise constant controls $(w, a, b)= (w(t), a(t), b(t)) \in \mathcal{S}_1$ with at most $L=3N$ switches, such that the flow map $\Phi_T$ generated by \eqref{eq:nodep1} satisfies  
    \begin{equation}\label{exactsimul}
    \Phi_T(x_i; w, a, b)=y_i,\quad\forall i\in[N].\end{equation}
\end{theorem}

The representation of a finite dataset, as considered thus far, is naturally related to the universal approximation problem. In the present framework, the latter is understood as the ability to approximate any continuous function on a compact set by means of the flow map associated with the NODE. More precisely, the following approximation result holds in $L^2$:

\begin{theorem}[Universal approximation theorem, $P=1$, \cite{domenecNODES}]\label{thm:2}
    Let the dimension $d \ge 2$,\( T > 0 \) be fixed, and consider a bounded set \( \Omega \subset \R^d \). Then, for any \( f \in L^2(\Omega; \R^d) \) and \( \varepsilon > 0 \), there exist piecewise constant controls $(w, a, b)= (w(t), a(t), b(t)) \in \mathcal{S}_1$ with a finite number of discontinuities, such that the flow map \( \Phi_T \) generated by \eqref{eq:nodep1} satisfies
    \begin{equation}
    \| \Phi_T - f \|_{L^2(\Omega)} < \varepsilon.
    \end{equation}
\end{theorem}

The proof combines the simultaneous controllability property of NODEs with approximation arguments based on simple functions defined on hyperrectangles, together with the ability of NODEs to compress and expand the regions where these simple functions take constant values.

In \cite{alvlop2024interplay}, we extend the analysis to any width $P \geq 1$ and show that the same result holds, with the number of needed switches $L$ decreasing as the width $P$ increases. Our findings reveal a trade-off between these two parameters: depth $L$, represented by the number of switches, and width $P$, as captured in the following result:

\begin{theorem}[Depth versus width, \cite{alvlop2024interplay}] \label{thm:pvsL}
    Let $d\ge 2$, \( N \geq 1 \) and \( T > 0 \) be fixed. Given \( P \geq 1 \) and \( \{ x_i, y_i \}_{i \in [N]} \subset \R^d \) with \( x_i \neq x_j \) and \( y_i \neq y_j \) for \( i \neq j \),  there exist piecewise constant controls $(W, A, b)= (W(t), A(t), b(t)) \in \mathcal{S}_P$ with \( L = 2\left\lceil N/P \right\rceil - 1 \) switches, such that the flow map \( \Phi_T \) generated by \eqref{eq:node} satisfies \eqref{exactsimul}.
\end{theorem}

The proof of this result is based on a two-step control method, applied inductively. The width  \( P \) enables parallelizing the movements in the inductive proof developed for  $P=1$. As \( P \) increases, the number of switches \( L \) decreases at the same rate, suggesting that both parameters play complementary roles in our control method. Remarkably, the overall complexity of the NODE, quantified by the product of the number of switches ($L$) and the number of neurons ($P$), remains invariant.

However, even when  \( P \geq N \), at least one switch (\( L = 1 \)) is required to transition between the two steps.  Thus, the above result does not yield autonomous dynamical systems, even when the number of neurons $P$ exceeds the cardinality of the dataset $N$. Nevertheless, by analogy with the results in the previous section on shallow NNs, one might expect similar outcomes to be achievable with autonomous (time-independent) NN vector fields.

We next study the autonomous NODE model, corresponding  to  \( L = 0 \) switches:
\begin{equation}\label{nodeauto}
    \dot{ x}(t) = W\vecsigma(Ax(t)+b),\hspace{1cm}t\in(0,T),
\end{equation}
in which, now, $(W, A, b) \in \R^{d\times P} \times \R^{P\times d} \times \R^P$ are time-independent.

We now relax the objective and prove the approximate simultaneous controllability property.
\begin{theorem}[Approximate simultaneous control for autonomous NODEs, \cite{alvlop2024interplay}]\label{autonomous}
    Let $d\ge 2$, $N\geq1$ and $T>0$ be fixed. Given  \(\{ (x_i, y_i)\}_{i\in[N]}\subset\mathbb{R}^{d}\times \mathbb{R}^{d}\) with $x_i\neq x_j$, for all $P\geq1$ there exists a constant control $(W,A,b)\in\mathbb{R}^{d\times P}\times\mathbb{R}^{P\times d}\times\mathbb{R}^P$ such that the autonomous NODE flow map $\Phi_T$  generated by \eqref{nodeauto} satisfies
    \begin{equation}\label{approxsimul}
        \sup_{i\in[N]} |y_i-\Phi_T(x_i)|\leq C\,\frac{\log_2(\kappa)}{\kappa^{1/d}},
    \end{equation}
    where $\kappa=(d+2)dP$ is the number of parameters in the model, and $C>0$ is a constant independent of $P$.  
\end{theorem}
Obviously, given $\varepsilon >0$, taking the width $P$, and therefore also $\kappa$, large enough, the right-hand side in \eqref{approxsimul} can be made smaller than $\varepsilon$, thus ensuring approximate simultaneous controllability.

\begin{figure}[h]
    \begin{minipage}{0.32\textwidth}
        \centering
        \includegraphics[width=\linewidth]{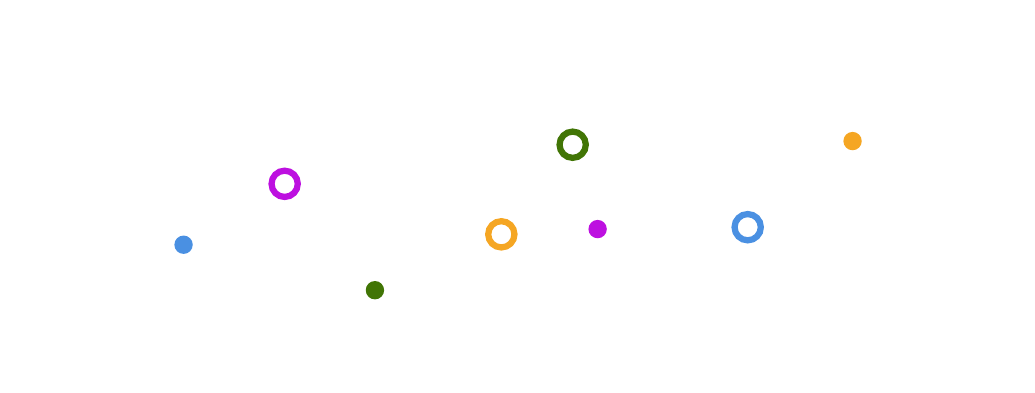}
    \end{minipage}
    \begin{minipage}{0.32\textwidth}
        \centering
        \includegraphics[width=\linewidth]{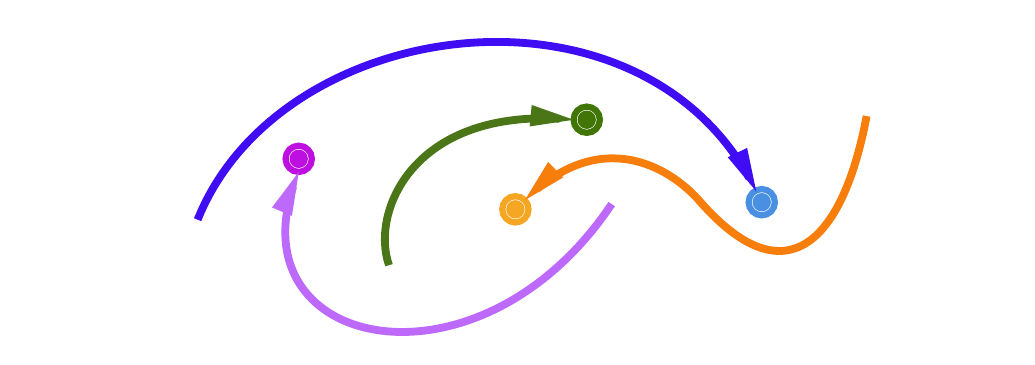}
    \end{minipage}
    \begin{minipage}{0.32\textwidth}
        \centering
        \includegraphics[width=\linewidth]{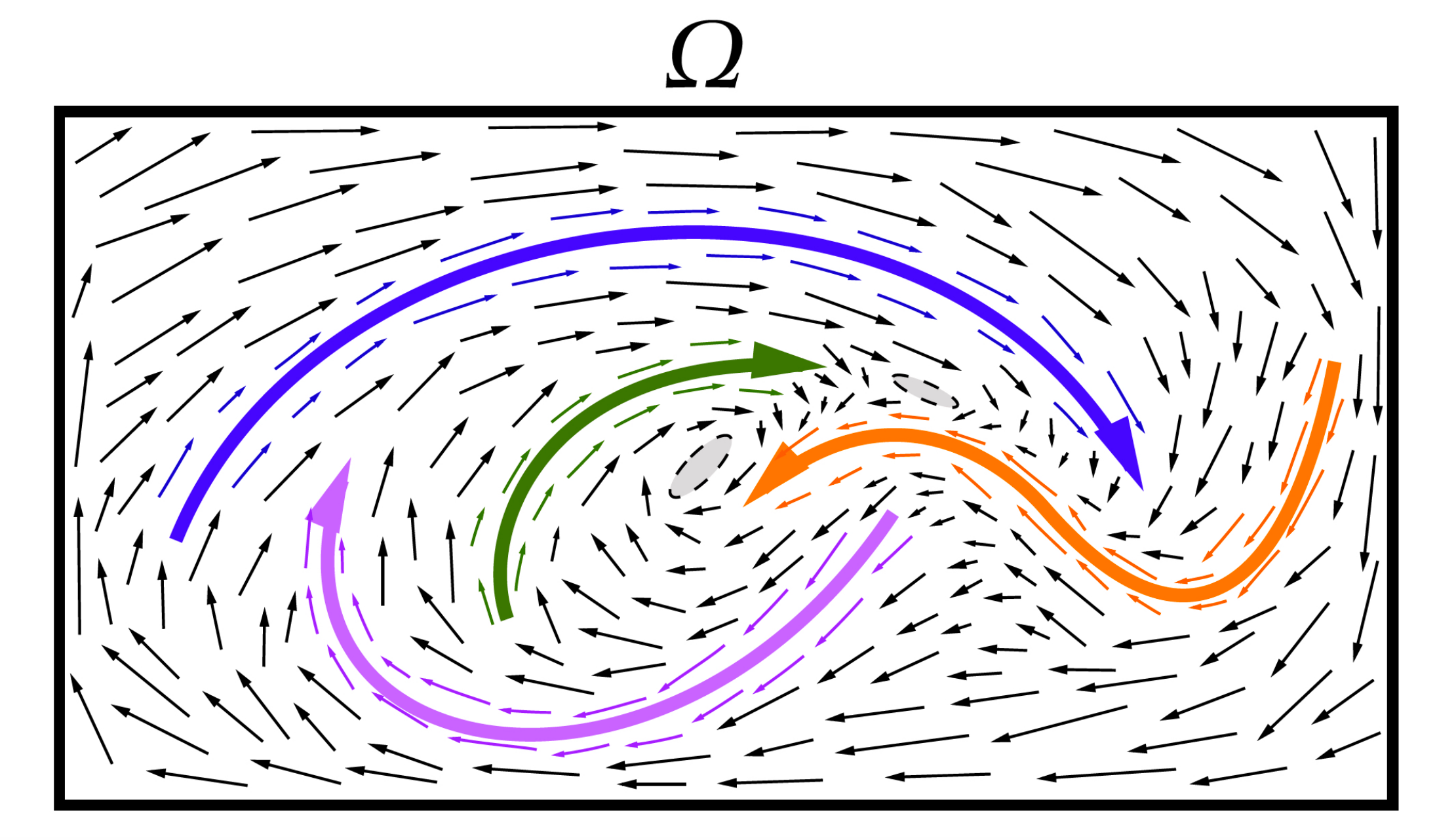}
    \end{minipage}
    \caption{From left to right: Construction of a vector field whose integral curves interpolate the dataset, defined in a compact domain $\Omega$ containing all the curves.}
    \label{fig:tubular}
\end{figure}

The proof of this result proceeds in two steps. First, we construct a smooth and globally Lipschitz time-independent vector field whose integral curves steer each input point  \( x_i\) to its corresponding target \( y_i\) within a fixed time \(T\) (assuming the dataset is consistent). The construction of this field is described in Figure \ref{fig:tubular}. In the second step, we apply the UAT to approximate the vector field by one generated through a shallow NN \cite{bach2017,devore_hanin_petrova_2021}. The final output is an approximate simultaneous control result, rather than an exact one, because in the second step, we rely on approximating the vector field.

Two considerations motivate autonomous NODEs. On the one hand, as previously noted and shown in the preceding section, this objective can be met using shallow NNs, which are static models. This suggests that temporal variation in the model is not necessary to achieve the required representational capabilities. On the other hand, it relates to the turnpike principle, whose origins can be traced back to John von Neumann, which ensures that optimal control strategies remain nearly constant over long time periods. We refer to \cite{geshkovski2022turnpike}, where this principle is applied to designing simpler and more stable architectures for deep ResNets.

To conclude this section, we present briefly some other related results:
\begin{enumerate}
      \item \emph{High dimensions}: If \(d > N\),  the number of switches can be reduced to   \(L=2\left\lceil N/P\right\rceil - 2\), \cite{alvlop2024interplay}.
    
    \item \emph{Probabilistic estimates}: We can also estimate the probability that the points will be arranged in spatial configurations that facilitate their autonomous control. For instance, if \( x_i\) and \( y_i\) are independently sampled from  \(U([0,1]^d)\) for all \(i\in[N]\), when $d\to+\infty$ for fixed $N$, with high probability,  autonomous NODEs suffice for exact representation,  \cite{alvlop2024interplay}.
    
    \item \emph{Impact of clustering}: The inductive constructions of controls enabling complete classification, fundamentally based on classifying each data point individually, can be substantially improved by exploiting the clustering properties of the dataset. In simplified terms, data clusters that can be separated by hyperplanes within a polytope are classified simultaneously through piecewise constant vector fields. Since randomly chosen datasets generally exhibit clustering,  the methods described above, when combined with a probabilistic analysis of dataset clustering, can reduce the number of switches and active neurons required in the NODE for classification. In \cite{alvlop2023optimized}, we pursue this direction by introducing new algorithms that incorporate the spatial structure of the data distribution to minimize the number of necessary parameters.
    \end{enumerate}

\subsection{The impact on training.}\label{impacttraining}
In the preceding subsections, we thoroughly examined the problem of data representation and classification through the lens of the simultaneous control properties of NODEs. The focus was on developing explicit constructions to identify parameter values that achieve the desired objectives. These constructions yield specific parameter configurations whose norms and number of switching instances can be quantitatively estimated.

This approach stands in contrast to the one commonly used in practice, where NNs are trained computationally by minimizing a cost functional based on empirical risk. Typically, the parameters obtained through training differ substantially from those derived via the geometric and inductive constructions presented earlier.

This duality between explicitly constructing parameters to demonstrate existence (with quantitative guarantees) and obtaining them computationally through optimization is not unique to NODEs for supervised learning. It also arises in classical control problems for both ODEs and PDEs. In such settings, the first approach emphasizes theoretical guarantees, including the existence of controls and the design of control strategies, while the second addresses the practical task of numerically computing optimal controls. Insights gained from the theoretical approach often provide valuable bounds and guidance for the numerical one. This interplay was, for example, analyzed in the context of approximate controllability of the heat equation in \cite{fernandez2000cost}.

The same applies in this context. To illustrate this important duality, let us go back to the simplest NODE in \eqref{eq:nodep1} endowed with one neuron, $P=1$. Fix a time horizon $T$ and a consistent dataset  \(\{ (x_i, y_i)\}_{i\in[N]}\subset\mathbb{R}^{d}\times \mathbb{R}^{d} \). We aim to build controls $(w(t), a(t), b(t))$ such that solutions of \eqref{eq:nodep1} map  $x_i$ to $y_i$ for each $i\in [N]$ in time $t=T$. We have shown that this can be done by constructive methods that allow us to estimate the complexity of the needed controls, in terms of the number of switching instances. Assume that the norms of the constructed controls are tracked carefully, as our methods do, to get, for some $K>0$, a bound of the form
\begin{equation}\label{explicitbound}
|| (a, b, w) ||^2_{BV(0, T; \R^{2d+1})} \le K.
\end{equation}

In practice, these constructions are typically disregarded, and training is instead carried out by minimizing a cost functional of the form
\begin{equation}\label{minimizationtask}
\min_{(a, b, w) \in BV(0, T; \R^{2d+1})} J_\alpha(a, b, w),
\end{equation}
\begin{equation} \label{thefunctional}
J_\alpha(a, b, w) = \alpha || (a, b, w) ||^2_{BV(0, T; \R^{2d+1})} + \sum_{i\in [N]} |\Phi_T(x_i)-y_i|^2. \end{equation}
Here $\alpha>0$ plays the role of a regularization or penalization parameter. The existence of global minimizers of $J_\alpha$  in \eqref{thefunctional} is easy to prove by the Direct Method of the Calculus of Variations \cite{brezis2011functional}. Note that the penalization term has been chosen in the $BV(0, T; \R^{2d+1})$-norm to avoid compactness issues and to assure that the piecewise constant controls built above are admissible.

The bound \eqref{explicitbound} is immediately useful. Indeed, whatever $\alpha >0$ is, any global minimizer of $J_\alpha$, that we denote by $(a^*, b^*, w^*) $, will necessarily satisfy the bound
\begin{equation}
||(a^*, b^*, w^*)||^2_{BV(0, T; \R^{2d+1})} \le K.
\end{equation}
This is simply because 
 \begin{equation}
 ||(a^*, b^*, w^*)||^2_{BV(0, T; \R^{2d+1})}
\le 
\frac{1}{\alpha} J_\alpha(a^*, b^*, w^*)  \le \frac{1}{\alpha} J_\alpha(a, b, w) \le  K,
\end{equation}
$(a, b, w)$ being the controls and parameter values we have built, for which, due to the simultaneous control property, the empirical risk  $\sum_{i\in [N]} |\Phi_T(x_i)-y_i|^2$ vanishes and the bound \eqref{explicitbound} holds.
 
This yields a priori insight for training via numerical optimization, a fact that
 applies not only to NODE architectures, but for all models we might employ and, in particular, to those discussed in these notes.
 
\subsection{Control of measures.}
Generative modeling and generalization are among the most impactful practical applications of NODEs. The goal is to learn the underlying data distribution of a given random variable, typically through sampling from a learned distribution and supervised learning, to generate new synthetic data or accurately assign labels to unseen samples.

 Normalizing flows \cite{papamakarios_normalizing_2021} are among the most effective methodologies for this purpose. In these models, a transformation is learned that maps a simple probability distribution, such as a uniform or standard Gaussian, into a complex target distribution, known only through a training dataset. Synthetic data can then be generated by sampling points from the simple distribution and applying the learned transformation.

NODEs play a key role in this paradigm, offering a new framework in which the flow is described as a continuous-time dynamic. To achieve this, the continuous model is reformulated in terms of transport equations, leveraging the classical relationship between \eqref{eq:node}, seen as the ODE of characteristics, and the corresponding hyperbolic transport PDE or continuity equation, i.e., the neural transport equation:	
\begin{equation}\label{eq:neuraltransport}
    \left\{
    \begin{array}{l}
        \partial_t \rho + \operatorname{div}_x(w(t) \sigma(\langle a(t), x\rangle + b(t))\rho) = 0,\\			
        \rho(0)=\rho_0,
    \end{array}\right.
\end{equation}
describing the evolution of a probability density $\rho_0\in L^1(\R^d)$.

We reformulate the problem from a controllability perspective: can we transform one probability measure into any other (possibly in an approximate manner, up to an arbitrarily small error) by selecting the appropriate controls $(w(t), a(t), b(t))$ for equation \eqref{eq:neuraltransport}? 

The error metric is crucial. The following result in  \cite{DomenecNormalisingFlows} guarantees approximate controllability in $L^1$.

\begin{theorem}[$L^1$-approximate control of neural transport, \cite{DomenecNormalisingFlows}]\label{thm:approxtv}
   Assume that $d\ge 2$.
    Given two probability densities \( \rho_0, \rho_T \in L^1(\R^d) \), for any \( T > 0 \) and \( \varepsilon > 0 \), there exist piecewise constant controls \( (w, a, b) \in \mathcal{S}_1 \) with a finite number of discontinuities, such that the solution of \eqref{eq:neuraltransport} with $P=1$ satisfies:
    \[
    \|\rho(T) - \rho_T\|_{L^1(\R^d)} < \varepsilon.
    \]
\end{theorem}

 Alternatively,  the question can be posed in terms of the Wasserstein distance, thereby linking this framework to optimal transport theory,  \cite{DomenecNormalisingFlows}. Under suitable conditions on \( \rho_0 \) and \( \rho_T \), Theorem \ref{thm:approxtv} was established in \cite{alvarez2025constructive} for the standard and commonly used Kullback--Leibler divergence. Extensions to arbitrary widths \(P \geq 1\) are studied in \cite{alvlop2024interplay}.

\subsection{Conclusions.}
This section reframed supervised learning within a control-theoretic framework, interpreting deep architectures as dynamical systems: ResNets as discrete dynamics and Neural ODEs as their continuous counterparts. Under this perspective, representation and classification become problems of ensemble controllability. This approach provides:
\begin{itemize}
    \item Constructive procedures with quantitative bounds for steering datasets,
    \item A clear depth-width trade-off captured by the balance between switches and neurons,
    \item Insight into autonomous flows, which play a special role in approximation,
    \item Complexity reduction via clustering, exploiting data structure,
    \item A natural extension from points to measures, linking NODEs to normalizing flows, transport metrics, and generative modeling.
\end{itemize}

Moreover, the constructive bounds obtained in this framework yield a priori estimates that directly inform practical training through regularized empirical risk minimization. Altogether, the control viewpoint supplies both conceptual clarity and useful design principles for deep learning, providing the foundation for the subsequent analysis of transformer architectures.

\section{Self-attention as a clustering mechanism in transformers.}\label{LLM}
We have analyzed the representation capacity of several architectures, with particular emphasis on NODEs. We demonstrated that the complexity of NODEs required to perform representation tasks can be quantified by the number of control switches, which constitutes an indicator that can, in turn, be estimated based on the size of the dataset to be represented. Furthermore, we showed that these estimates can be sharpened by exploiting the clustering structure of the data. This naturally raises the question of whether modern architectures, such as transformers, inherently facilitate such clustering.

We show that self-attention in transformers (on which the most recent and capable large language models (LLMs) are based) enhances clustering, thereby reducing the complexity of the NODEs needed for representation. This provides a partial, yet mathematically meaningful, interpretation of the utility of self-attention: it enables a reduction in architectural complexity for deep NNs in supervised learning tasks.

\subsection{Dynamics of self-attention layers and their asymptotics.}
Transformers demonstrate superior performance in supervised learning tasks, largely due to their ability to capture context; that is, the relationships between words within a sentence. To achieve this, transformers are trained on datasets consisting of sequences of words, such as sentences or paragraphs. Formally, the training data is given by $\{(Z^j, Y^j)\}_{j\in[N]}$, where each input sequence $Z^j \in (\mathbb{R}^d)^n$ and output sequence $Y^j \in (\mathbb{R}^d)^m$ represent $n$ and $m$ words encoded as vectors in a $d$-dimensional Euclidean space. Often in applications, $m \leq n$, so that the length of the target sequences $Y^j$ is smaller than that of the input sequences $Z^j$.

To effectively capture contextual relationships within sequences, the transformer architecture is often employed in the state-of-the-art  LLMs. Transformers, which can be viewed as an extension of ResNets, incorporate self-attention layers that exploit the sequential structure of the data. Heuristically, these layers model the ``context'' of each sequence by dynamically weighting and mixing its components based on pairwise similarity. 

Motivated by their empirical success, a rigorous theoretical framework for understanding the role and mechanisms of self-attention is beginning to emerge. Early work in this direction includes \cite{lu2019understanding,sanderSinkformers2022}, which interprets the elements of an input sequence as interacting particles, where the interaction is mediated by a kernel derived from the self-attention mechanism. This particle-based viewpoint is further developed in \cite{geshkovski2023emergence,burger2025analysis,geshkovski2025mathematical}, where it is used to establish asymptotic clustering results for attention-only transformers with shared weights.

Our work \cite{alcalde2025clustering} adds to this growing literature by providing precise mathematical results that explain the role of self-attention as a clustering mechanism, within a simplified yet expressive class of attention-only transformers which we refer to as \emph{hardmax self-attention dynamics}.

Such dynamics are parameterized by a symmetric positive definite matrix $A \in \mathbb{R}^{d\times d}$ and a scalar parameter $\alpha > 0$. They act on components $\zz{i}{} \in \R^d$ of a sequence $Z=(\zz{1}{},\ldots,\zz{n}{})\in (\mathbb{R}^d)^n$, called \textit{tokens} in the ML literature. Given initial token values $Z(0)=(\zz{1}{}(0),\ldots,\zz{n}{}(0))$, and denoting by $Z(k)$  the sequence of tokens in layer $k$, they  evolve according to the following discrete-time dynamics:\begin{subequations}\label{eq:transformer}
\begin{equation}
\label{eq:transformer_a}
\zz{i}{}(k+1)=\zz{i}{}(k)+\frac{\alpha}{1+\alpha}\,\frac{1}{|\mathcal{C}_{i}(k)|}\sum_{\ell\in \mathcal{C}_{i}(k)}\Big(\zz{\ell}{}(k)-\zz{i}{}(k)\Big),
\end{equation}
\begin{equation}
\label{eq:transformer_b}
\mathcal{C}_{i}(k)= 
    \left\{ 
    j\in [n]\,\, : \,\,
    \big\langle \zz{j}{}(k) , A \zz{i}{}(k) \big\rangle
    =
    \max_{\ell\in[n]}
    \big\langle \zz{\ell}{}(k) , A \zz{i}{}(k) \big\rangle
    \right\},
\end{equation}
\end{subequations}
where $|\mathcal{C}_{i}(k)|$ denotes the cardinality of the index set $\mathcal{C}_{i}(k)$.

\begin{figure}
    \centering
    \begin{subfigure}{.5\textwidth}
      \centering
      \includegraphics{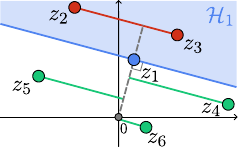}
      \caption{}
    \end{subfigure}%
    \begin{subfigure}{.5\textwidth}
      \centering
      \includegraphics{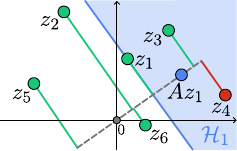}
      \caption{}
    \end{subfigure}
    \caption{Geometric interpretation of \eqref{eq:transformer_b} for $i=1$ with (a) $A=I$ and (b) $A = \left(\begin{smallmatrix} 2 & 1 \\ 1 & 1 
    \end{smallmatrix}\right)$. In (a), tokens $\zz{2}{}$ and $\zz{3}{}$ have the largest orthogonal projection on the direction of $A \zz{1}{} = \zz{1}{}$, so $\mathcal{C}_1(Z) = \{ 2,3 \}$. In (b), token $\zz{4}{}$ has the largest projection on the direction of $A\zz{1}{}$, so $\mathcal{C}_1(Z) = \{ 4 \}$. In both cases, tokens attracting $z_1$ can only lie on the closed half-space $\mathcal{H}_1 = \{z:\;\langle A \zz{1}{}, \zz{}{} - \zz{1}{} \rangle \geq 0\}$ (blue shading).}
    \label{fig:geom_interp}
\end{figure}

Viewed as a discrete-time dynamical system describing the evolution of tokens, \eqref{eq:transformer} has a simple geometric interpretation: token $\zz{i}{}$ is attracted to the average of the tokens with the largest orthogonal projection in the direction of $A \zz{i}{}$ (cf. Figure \ref{fig:geom_interp}),  $\alpha$ being  a \textit{step-size} parameter regulating the intensity of the attraction.

Given the depth of modern transformers, we analyze the asymptotic behavior of tokens that evolve according to \eqref{eq:transformer}. Fixing $\alpha$ and $A$, we prove that as $k\to \infty$ tokens converge to a clustered equilibrium constituted either by special tokens, which we call \emph{leaders}, or particular convex combinations thereof. More precisely,
token $\zz{i}{}$ is a \emph{leader} if there exists $k_i \in \N$ such that $\mathcal{C}_{i}(k) = \{ i d\}$ for all $k \geq k_i$. We denote as $\mathcal{L}(k) = \{ z_i(k) \in Z(k) : \; \mathcal{C}_{i}(k) = \{ i\} \}$ the set of leaders at layer $k$. 

The following asymptotic behavior result was proved in \cite{alcalde2025clustering}:

\begin{theorem}[Asymptotic clustering, \cite{alcalde2025clustering}]
\label{thm:emergenceClusters}
Assume the initial tokens $Z(0)=(\zz{1}{}(0),\ldots,\zz{n}{}(0)) \in (\R^d)^n$ to be nonzero and distinct, the matrix $A\in \R^{d\times d}$ in \eqref{eq:transformer_b} to be symmetric and positive definite, and $\alpha >0$. 

Then, as $k \to \infty$, $Z(k)$ converges to a clustered equilibrium configuration $Z^*$ with the property that each of its tokens is a leader or a projection (with respect to the norm associated with $A$) of the origin onto a face of $\mathcal{K}$, the convex polytope with the leaders as vertices. 

Furthermore, the set of leaders stabilizes after a finite number of layers, i.e., there exists $k_0$ such that $\mathcal{L}(k) \equiv \mathcal{L}(k_0)$ for all $k \geq k_0$.
\end{theorem}

The proof uses a geometric analysis of the discrete dynamics induced by the self-attention model under consideration. The asymptotic stability arises from the balance between two opposing geometric effects: on the one hand, the norm of each token increases, while on the other, the convex envelope of the tokens contracts in the sense of set inclusion. The only possible reconciliation of these two competing forces is convergence to a clustered equilibrium of the type described above.

    Several remarks about the previous result are in order.
    \begin{itemize}
        \item Although the set of leaders stabilizes after finitely many layers (as guaranteed by Theorem \ref{thm:emergenceClusters}), convergence of all tokens to their final positions is only asymptotic. The number of layers required to approach equilibrium is highly sensitive to the initial configuration, and no uniform bound is available.  
        
                \item Leaders emerge during the evolution of the dynamics, although in special cases, an initial token value may remain a leader throughout. The number of distinct leaders that ultimately emerge depends on the initial configuration of tokens, and is always at most $n$. 
        \item  At a clustered equilibrium, each token is either a leader or a point on a face of the polytope spanned by the leaders. Consequently, the effective asymptotic complexity of the system is determined by the number of distinct leaders that emerge.
    \end{itemize}

We refer to \cite{alcalde2025clustering} for a computational experiment on a supervised sentiment analysis task which exploits the asymptotic simplification properties guaranteed by the result above.

\subsection{Exact interpolation of sequences with transformers.}
The theoretical clustering results discussed above can also serve to enhance the representational capacity of neural network architectures. Building on this idea, in \cite{alcalde2025exact} we investigate the interpolation capabilities of a broader class of transformers, consisting of alternating ResNet and self-attention layers, which we interpret as a discrete dynamical system.

Since we identify sequences with finite sets of elements in $\R^d$, two sequences are regarded as equal if they differ only by a permutation of their elements, a property known as permutation invariance. Moreover, we assume the following properties for the dataset of sequences: i) the input and target sequences $\{(Z^j, Y^j)\}_{j\in[N]}$ satisfy that the sequences $Z^1,\ldots, Z^N$ are pairwise distinct, and ii) tokens within each input and output sequence are pairwise distinct. Let $\mathrm{R}$ denote a ``readout'' map, which we fix as the set-to-set transformation that removes repeated elements from the input set. For instance, for $z_1 \neq z_2$, then $\mathrm{R}(\{z_1, z_1, z_2\}) = \{ z_1, z_2\}$.

We pose the following \textit{exact interpolation} problem: find a transformer $\T : (\R^d)^n \to (\R^d)^n$ such that $(\mathrm{R} \circ \T)(Z^j) = Y^j$ for all $j \in [N]$. 

To illustrate this framework, we now present two classical examples of exact interpolation.
The first is \emph{next-token prediction}, where each sequence $Z^j$ corresponds to an incomplete sentence and the target $Y^{j}=\{y^{j}\}$ is a single token encoding the missing word, i.e., $m=1$.
The second is \emph{sequence-to-sequence interpolation}, where the target $Y^j = \{ y_1^{j}, \dots, y_{m}^j \}$ represents the possible set of words completing the sentence.

In \cite{alcalde2025exact}, we show that the exact interpolation problem can be solved with transformers alternating ResNets (with the ReLU activation function) and self-attention steps. Furthermore, we provide explicit complexity bounds, establishing that transformers need a total of $\mathcal{O}(N\, m)$ ResNet and self-attention layers and $\mathcal{O}(d \, N \, m)$ nonzero parameters to achieve exact interpolation.
One of the remarkable strengths of this transformer architecture is that its parameter complexity is independent of the input sequence length $n$, but depends on the output sequence length $m$. This stands in clear contrast to pure ResNets or NODEs, where, as we have seen, the number of parameters scales linearly with the dimensions of the dataset, which in the present setting is $\mathcal{O}(d\, n\, N)$.

Our results, initially derived for transformers with a hardmax self-attention formulation, also apply to the standard softmax self-attention. The latter is defined with a temperature parameter $\tau > 0$ as
\begin{equation}
    \label{eq:softmaxFormulation}
   \pi_{i\ell}^\tau (Z) = \frac{\exp{\left(\frac{1}{\tau}\langle A z_i, z_\ell \rangle\right)}}{\sum_{k=1}^{n} \exp{\left(\frac{1}{\tau}\langle A z_i,  z_k \rangle\right)}}.
\end{equation}
As $\tau \to 0$, the singular limit leads to the  hardmax self-attention formulation:
\begin{equation}
    \label{eq:hardmaxFormulation}
   \pi_{i \ell }^0 (Z) =
    \begin{cases}
        \frac{1}{ | \mathcal{C}_i(Z) |} &\text{if } \ell \in \mathcal{C}_i(Z),\\
        \;\, 0 &\text{otherwise}
    \end{cases}
\end{equation}
where $\mathcal{C}_i(Z)$ is defined as in \eqref{eq:transformer_b} above.

Our results are proved constructively with a simultaneous control strategy, characterizing explicitly the parameters of each layer. Such a strategy leverages self-attention layers in two key ways: first, by allowing sequences to become pairwise disjoint, and second, by clustering the $n$ tokens of the input sequence to the $m$ tokens of the target sequence. This effectively removes the dependency of the total number of parameters required for exact interpolation on the input sequence length $n$, typically large. The number of required parameters is then governed by the output sequence length $m$.

We first prove the result for transformers with the hardmax self-attention formulation. The extension to the softmax case requires one additional step: the inclusion of nonresidual ResNet layers to guarantee exact representation, since the intrinsic regularizing effect of the softmax makes this otherwise difficult to achieve.

Our exact interpolation results also provide explicit training bounds in the spirit of Subsection \ref{impacttraining}.

\subsection{Conclusions.} Our analysis establishes self-attention as a natural clustering mechanism in transformer architectures. By interpreting attention dynamics through a geometric, particle-based framework, we showed that tokens asymptotically converge to clustered equilibria characterized by a finite set of leaders. This clustering reduces the effective complexity of the representation problem and enables the design of transformers that achieve exact sequence interpolation with explicit parameter bounds independent of input sequence length, but dependent on target sequence length. Altogether, these results provide a rigorous explanation of how self-attention improves efficiency in modern architectures, linking theoretical principles to the empirical success of transformers.

\section{Related topics and perspectives.}
We outline related topics, each of which constitutes a subject worthy of investigation, rich in challenging open problems. Our topical selection is necessarily limited and reflects themes influenced by our recent work. These, together with many others, form a fascinating landscape to be explored with the overarching goal of merging classical methodologies in computational and applied mathematics with those inspired by ML. For further details, we refer to \cite{zuazuaober}.

\begin{itemize} 
\item  \emph{Time Reversal in Diffusion Models for Generative AI. } 
The ideas and methods developed in Section \ref{SNN} can also be applied to the time inversion of the heat equation (see \cite{liu2025moments}), building on the idea of reducing its infinite-dimensional, highly irreversible dynamics to a finite-dimensional system that tracks only a finite number of solution moments. The time inversion of heat processes is also a key ingredient in modern generative AI techniques, a fascinating connection that we briefly highlight here.
 
 The classical Li--Yau inequality (\cite{liyau}) for positive solutions \(u \ge 0\) of the \(d\)-dimensional heat equation,  
\[
u_t - \Delta u = 0 \quad \text{in } \mathbb{R}^d \times [0,\infty),
\]  
ensures that  
\[
\Delta \log(u) \;\ge\; -\frac{d}{2t}.
\]

Given a dataset \(\{x_i\}_{i=1}^I\) sampled from an unknown distribution, we define the empirical measure  and the corresponding parabolic solution 
$$
u_0(x) = \frac{1}{I}\sum_{i=1}^I  \delta(x - x_i), \quad
u(x,t) = \frac{1}{I} \sum_{i=1}^I \, G(x - x_i,t),
$$
$G$ being the Gaussian heat kernel in $\mathbb{R}^d$, i.e., $G(x, t)= (4\pi t)^{-d/2}\exp(-|x|^2/4t)$, which diffuses the information of the initial empirical measure throughout \(\mathbb{R}^d \times [0,\infty)\).  

Generative diffusion models aim to reverse this diffusion process to generate new samples from the same distribution. However, this backward evolution is severely ill-posed.   But this instability explains, partly,  their strong generative capacity.  
Our first contribution in this context is to show that the well-posedness of this inversion mechanism relies on Li--Yau's inequality above.

Observe that the backward heat equation can be rewritten as  
\[
u_t - \Delta u  = u_t + \Delta u - 2\Delta u
= u_t + \Delta u - 2 \, \text{div}\!\left(u \frac{\nabla u}{u}\right) 
= u_t + \Delta u - 2 \, \text{div}(u \nabla \log u) = 0,
\]  
and, introducing the {\it score function}, 
$
s(x,t) = \nabla \log(u),
$
it takes the form of a convection-diffusion or Fokker--Planck model:  
\[
u_t + \Delta u - 2 \, \text{div}\!\big(s(x,t) u\big) = 0,
\]  
which, unlike the original backward heat equation, is well-posed backward in time thanks to the Li--Yau inequality, which, in terms of the score function, reads:  
\[
\text{div}\, s(x,t) \;\ge\; -\frac{d}{2t}.
\]  
Indeed, this unilateral bound allows us to perform  the following energy estimate backward in time:  
\begin{equation*}
\frac{1}{2}\frac{d}{dt}\int u^2\,dx - \int |\nabla u|^2\,dx
= 2 \int \text{div}(u s(x,t)) \, u \, dx 
= \int u^2 \, \text{div}\, s(x,t)\, dx 
\;\ge\; -\frac{d}{2t} \int u^2 \, dx,
\end{equation*}  
yielding a priori estimates for backward solutions for all time \(t=\tau >0\). This estimate blows up as \(t \to 0^+\).  

The success of diffusion-based generative AI models relies partly on such well-posedness properties:  new samples of the unknown distribution are generated by realizing specific trajectories of the underlying backward stochastic differential equation. 

In this context, several hyperparameters may be tuned to regulate the generation capacity of the process: 
\begin{figure}[h]
    \centering
    \includegraphics[width=0.30\linewidth]{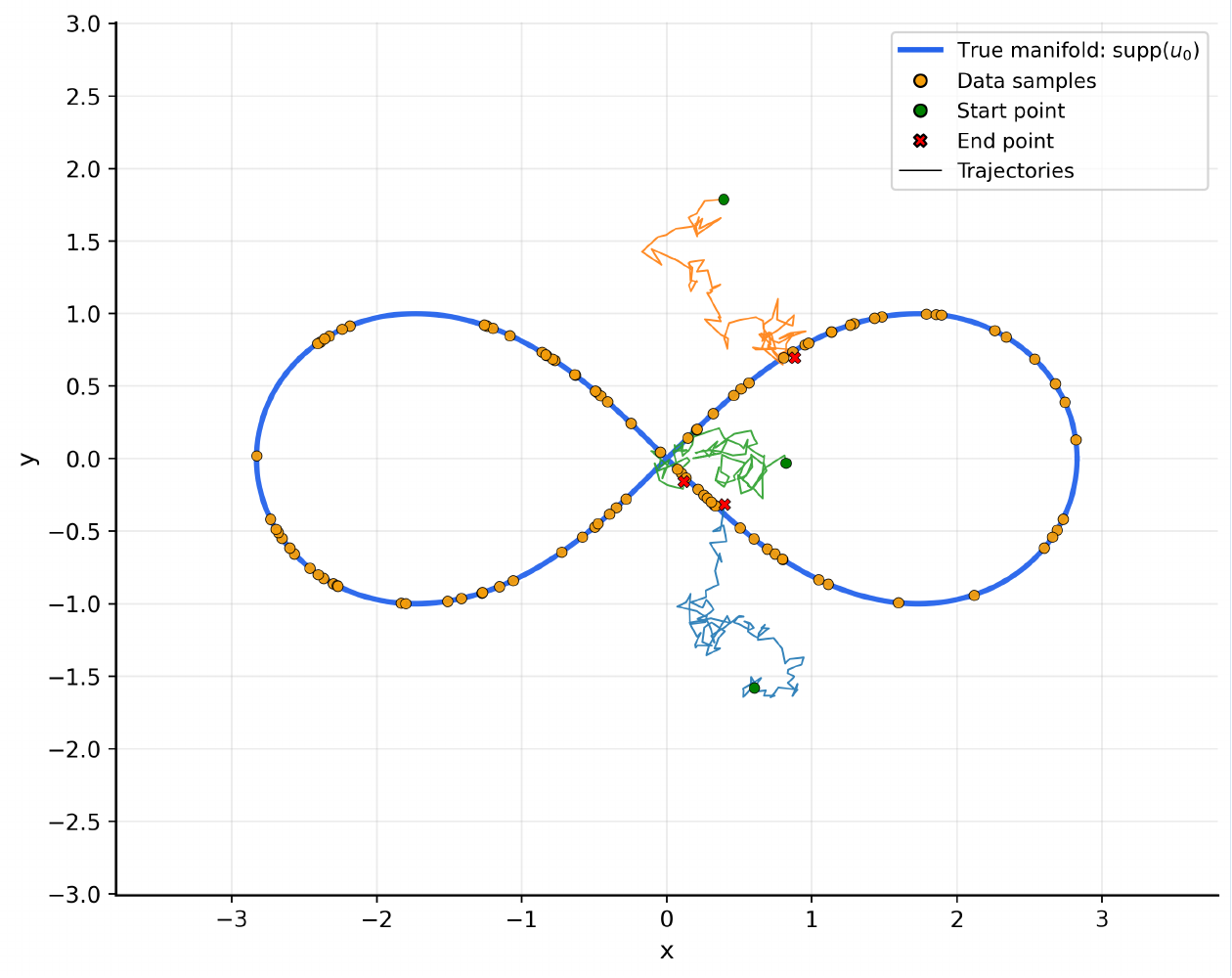} \hskip 1cm
    \includegraphics[width=0.30\linewidth]{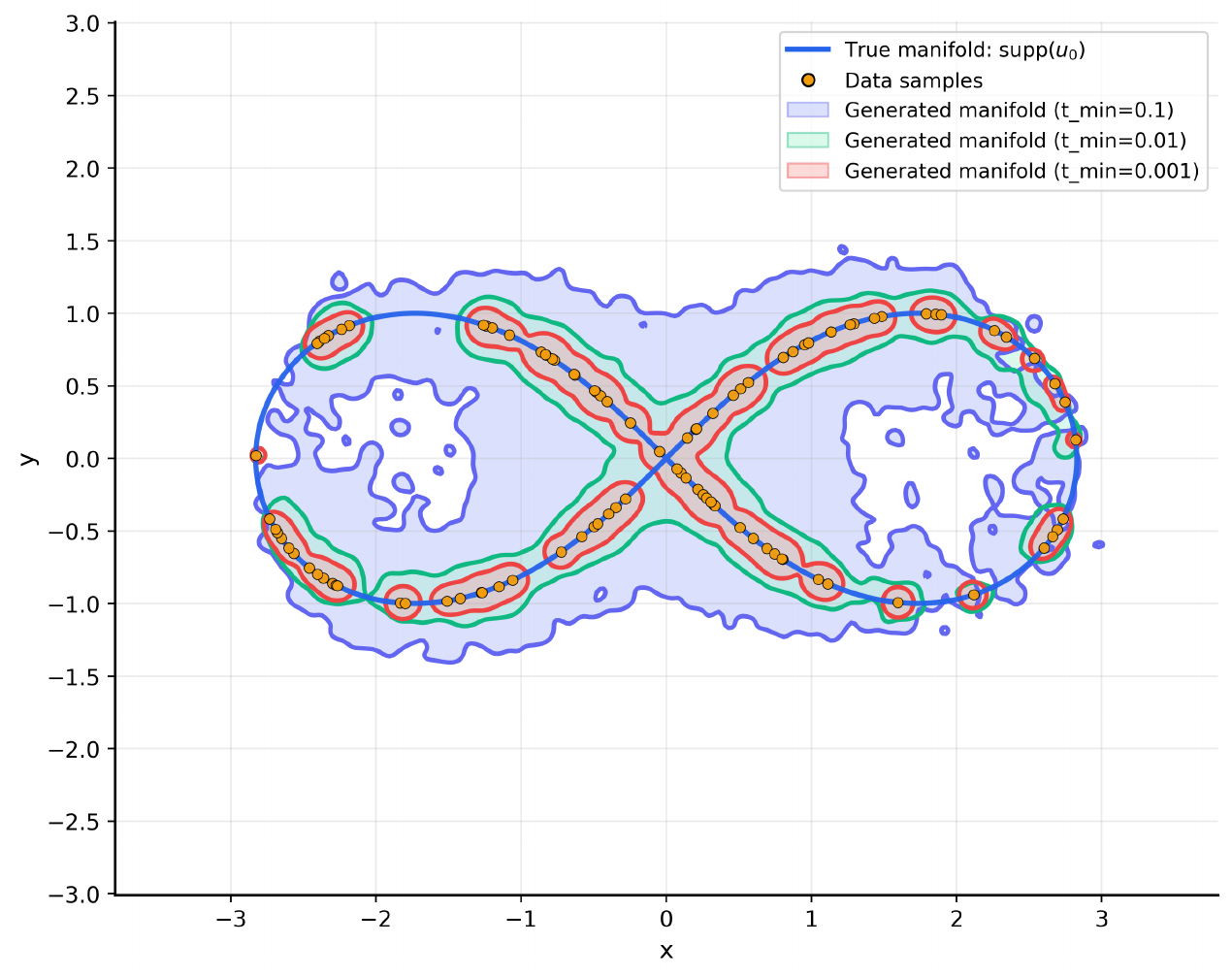}
    \caption{Left: Unknown data-manifold (blue lemniscate) and observed samples (orange). Three trajectories (with different initial points) of the backward stochastic differential equation driven by the score function determined by the Gaussian ansatz associated with the given samples. 
Right: Regions enclosing 10,000 generated points for different stopping times \( t_{\text{min}} \in \{0.1, 0.01, 0.001\}\).
}
    \label{fig:DM}
\end{figure}
\begin{itemize}
\item The stopping time \(t = \tau > 0\) serves as a design parameter: a smaller \(\tau\) produces samples closer to the initially available ones, while a larger \(\tau\) generates more diverse and distinct samples. A toy example in two dimensions illustrating this phenomenon is illustrated in Fig.~\ref{fig:DM}

\item   The diffusivity level \(\nu=\nu(t)\ge 0\) in the backward diffusion problem can also be tuned, since  
\[
- \Delta u = \nu(t) \Delta u - (1+\nu(t)) \, \text{div}(s(x,t) u).
\]  
\item
Finally, the exact Gaussian score function above can be approximated by an NN ansatz, yielding a surrogate \(\overline{s}(x,t)\). The estimates above highlight the importance of maintaining, for some finite \(M>0\), uniform bounds of the form  
$
\text{div}\, \overline{s}(x,t) \;\ge\; -M.
$
\end{itemize}

The connection between Li--Yau type parabolic inequalities and diffusion-based generative AI is both natural and profound, and will be further explored in \cite{liu2025generativeai}.

\item \emph{Multilayer perceptrons:} The dynamical control methods developed to explain and analyze supervised learning in ResNets and NODEs can also be extended to other NN architectures. For instance,  in \cite{Hernandez2024Classification}, we analyzed the multilayer perceptron (MLP) models of the form:
\begin{align}\label{eq:multilayer_perceptron}
    x^{k+1} = \vecsigma_{k+1}(A^k x^k + b^k), \quad k \in [L],
\end{align}
where
 $x^k\in \mathbb{R}^{d_k}$ is the state at layer $k\in[L]$; 
 $\{d_{k}\}_{k\in [L]}$ is a prescribed sequence of widths for each layer  $k\in[L]$; 
 $A^k \in \mathbb{R}^{d_{k+1} \times d_{k}}$ and $b^k \in \mathbb{R}^{d_{k+1}}$ denote the model parameters; and 
 $\vecsigma_{k+1}: \mathbb{R}^{d_{k+1}} \to \mathbb{R}^{d_{k+1}}$ is the component-wise ReLU activation function.

In \cite{Hernandez2024Classification}, we developed constructive methods that provide explicit parameter values, layer dimensions $d_k$, and depths required for exact interpolation, together with the corresponding bounds.
These results, in turn, are analogous to what we observed for other architectures, making it possible to derive explicit bounds on the parameters obtained through training in standard computational practice.

\item \emph{Federated Learning.}  Federated learning (FL) arises from a simple tension: the most valuable data for training modern models (mobile interactions, clinical records, financial transactions, etc.) are distributed across many owners and often too sensitive or costly to centralize. Regulations, institutional policies, bandwidth limits, and business concerns all discourage raw data sharing, yet organizations still want the accuracy benefits of training on diverse, large-scale datasets. FL resolves this by moving computation to the data: clients train locally and share only model updates, enabling collaborative learning without transferring raw records. This preserves privacy, reduces communication of large datasets, accommodates heterogeneous devices and non-IID data, and supports personalization, while unlocking cross-silo collaboration (e.g., between hospitals or banks). 

This topic is closely related to splitting, alternating, and domain decomposition methods in classical computational mathematics, and the know-how in these more mature areas can contribute to the development of improved FL methodologies. In \cite{wang2023approximatea}, we address the crucial issue of privacy and vulnerability to attacks; in \cite{song2024fedadmm}, we propose a self-adaptive version of FL that enhances performance while reducing computational cost; and in \cite{liu2024fedgame}, we provide a refined game-theoretical analysis of FL methodologies.

\item \emph{Model Predictive Control and Reinforcement Learning.}  The CT-ML interaction has intensified over the past decades, building on their complementary strengths. Two paradigms stand: Model Predictive Control (MPC) and Reinforcement Learning (RL).
MPC and RL share the goal of optimizing sequential decision-making through feedback, but differ in that MPC relies on explicit system models and real-time optimization, whereas RL learns policies directly from data, often without a model.

In \cite{veldman2022framework},  the convergence of MPC is analyzed when applied to a finite-dimensional Linear Quadratic Regulator (LQR) problem. The convergence of MPC relies on the fact that, over infinite horizons, optimal controls and solutions are characterized by the algebraic Riccati feedback operator, ensuring exponential decay as time grows. This property underlies the turnpike principle: as the horizon expands, optimal controls and trajectories converge to the infinite-horizon steady state.
The techniques in \cite{veldman2022framework} could be adapted to learning tasks and NN architectures, and may also contribute to the analytical foundations of RL, although such extensions would require substantial effort.
With this aim, \cite{veldman2024stability} further combines this analysis with the Random Batch Method (RBM), introduced in \cite{jin_random_2020} for particle systems. RBM reduces computational cost while ensuring convergence in expectation. Its integration with MPC is natural and was successfully applied in \cite{ko2021model} to a guidance-by-repulsion control problem for collective dynamics. It would be interesting to investigate the application of these methods in the RL setting, particularly in a data-driven framework.

\item \emph{NODE approximation of nonautonomous ODEs.}  The UAT can also be used to derive accurate approximations of ODEs and PDEs. In \cite{li2024universal}, we introduce the class of semi-autonomous neural differential equations (SA-NODE)  of the form
\begin{equation}\label{eqsanode}
 \dot{x}(t) = W  \vecsigma \big(A x(t) + \beta t + b\big).
\end{equation}
Here the state is vector valued, i.e., $x(t)\in \R^d$, and the coefficients $(W, A, \beta, b) \in \R^{d\times P} \times \R^{P\times d}\times \R^P \times \R^P$ are time-independent, $P\ge 1$ being the number of neurons involved. This ansatz enables the efficient approximation of the dynamics of general finite-dimensional nonautonomous systems.
The ansatz \eqref{eqsanode} is nonautonomous, though its coefficients are time-independent. In fact, the sole source of time dependence in the model is the vector-valued coefficient $\beta$, which multiplies the time variable in the argument of the activation function. This significantly diminishes the computational cost of training.
Moreover, in \cite{li2024universal}, motivated by the classical connection between transport equations and their characteristic ODEs, we address linear transport equations arising in fluid mixing.

\item \emph{Hybrid PDE and data-driven modeling methods.} To effectively hybridize classical computational techniques with data-driven ones, in \cite{lorenzothysez} we introduced the Hybrid--Cooperative Learning (HYCO) methodology. HYCO is a novel hybrid learning paradigm designed to combine physics-based and data-driven approaches while mitigating their respective limitations. Unlike well-established methods, such as physics-informed neural networks (PINNs) \cite{raissi2019physics}, which enforce physical and data constraints directly on a single synthetic model, HYCO takes a fundamentally different route. It trains two models in parallel: one grounded in physical principles, typically formulated as an ODE or PDE (the physical model), and another driven by data, using a NN ansatz (the synthetic model). Instead of merging the physics-based and machine-learning components into a single architecture, HYCO lets them interact, much like two experts exchanging insights before reaching consensus, in a game-theoretical setting. HYCO opens a new perspective and raises fundamental questions regarding convergence guarantees with rates that merit careful investigation.

\end{itemize}

\section*{Acknowledgments.}
The author thanks all team members and collaborators who generously contributed to the development of the research summarized here.  Our research has been funded by the European Research Council (ERC) under the European Union's Horizon 2030 research and innovation program (grant agreement no. 101096251-CoDeFeL), the Alexander von Humboldt--Professorship program, the ModConFlex Marie Curie Action, HORIZON-MSCA-2021-dN-01, the Transregio 154 Project of the DFG, AFOSR Proposal 24IOE027, and grants PID2020-112617GB-C22 and TED2021131390B-I00 of the AEI (Spain), and Madrid Government-UAM Agreement for the Excellence of the University Research Staff in the context of the V PRICIT.

\bibliographystyle{siamplain}
\bibliography{references}

@inproceedings{
zhang2017understanding,
title={Understanding deep learning requires rethinking generalization},
author={Chiyuan Zhang and Samy Bengio and Moritz Hardt and Benjamin Recht and Oriol Vinyals},
booktitle={International Conference on Learning Representations},
year={2017}}

@article{liu2025generativeai,
  title={A PDE Perspective on Generative Diffusion Models},
  author={Liu, Kang and Zuazua, Enrique},
 journal={arXiv preprint  arXiv.2511.05940},
  year={ 2025}
}

@article {ko2021model,
    AUTHOR = {Ko, Dongnam and Zuazua, Enrique},
     TITLE = {Model predictive control with random batch methods for a
              guiding problem},
   JOURNAL = {Math. Models Methods Appl. Sci.},
  FJOURNAL = {Mathematical Models and Methods in Applied Sciences},
    VOLUME = {31},
      YEAR = {2021},
    NUMBER = {8},
     PAGES = {1569--1592},
      ISSN = {0218-2025,1793-6314},
   MRCLASS = {93C15 (49M29 49N10 90C59 93B45 93B51)},
  MRNUMBER = {4307003}
}

@article {veldman2022framework,
    AUTHOR = {Veldman, D. W. M. and Zuazua, E.},
     TITLE = {A framework for randomized time-splitting in linear-quadratic
              optimal control},
   JOURNAL = {Numer. Math.},
  FJOURNAL = {Numerische Mathematik},
    VOLUME = {151},
      YEAR = {2022},
    NUMBER = {2},
     PAGES = {495--549},
      ISSN = {0029-599X,0945-3245},
   MRCLASS = {49N10 (37M05 49K45 65C99 65L05 65L20)},
  MRNUMBER = {4433122},
MRREVIEWER = {Iman\ Malmir}
}

@article {jin_random_2020,
    AUTHOR = {Jin, Shi and Li, Lei and Liu, Jian-Guo},
     TITLE = {Random batch methods ({RBM}) for interacting particle systems},
   JOURNAL = {J. Comput. Phys.},
  FJOURNAL = {Journal of Computational Physics},
    VOLUME = {400},
      YEAR = {2020},
     PAGES = {108877, 30},
      ISSN = {0021-9991,1090-2716},
   MRCLASS = {82C22 (65M75 82C31)},
  MRNUMBER = {4019084}
}

@article {veldman2024stability,
    AUTHOR = {Veldman, Dani\"el W. M. and Borkowski, Alexandra and Zuazua,
              Enrique},
     TITLE = {Stability and convergence of a randomized model predictive
              control strategy},
   JOURNAL = {IEEE Trans. Automat. Control},
  FJOURNAL = {Institute of Electrical and Electronics Engineers.
              Transactions on Automatic Control},
    VOLUME = {69},
      YEAR = {2024},
    NUMBER = {9},
     PAGES = {6253--6260},
      ISSN = {0018-9286,1558-2523},
   MRCLASS = {93B45 (49N10)},
  MRNUMBER = {4792609}
}

@article {weinan2017proposal,
    AUTHOR = {E, Weinan},
     TITLE = {A proposal on machine learning via dynamical systems},
   JOURNAL = {Commun. Math. Stat.},
  FJOURNAL = {Communications in Mathematics and Statistics},
    VOLUME = {5},
      YEAR = {2017},
    NUMBER = {1},
     PAGES = {1--11},
      ISSN = {2194-6701,2194-671X},
   MRCLASS = {49M20 (37N35 68T05 93A30)},
  MRNUMBER = {3627592}
}

@article {haber2017stable,
    AUTHOR = {Haber, Eldad and Ruthotto, Lars},
     TITLE = {Stable architectures for deep neural networks},
   JOURNAL = {Inverse Problems},
  FJOURNAL = {Inverse Problems. An International Journal on the Theory and
              Practice of Inverse Problems, Inverse Methods and Computerized
              Inversion of Data},
    VOLUME = {34},
      YEAR = {2018},
    NUMBER = {1},
     PAGES = {014004, 22},
      ISSN = {0266-5611,1361-6420},
   MRCLASS = {68T05},
  MRNUMBER = {3742361}
}

@article {li2018maximum,
    AUTHOR = {Li, Qianxiao and Chen, Long and Tai, Cheng and E, Weinan},
     TITLE = {Maximum principle based algorithms for deep learning},
   JOURNAL = {J. Mach. Learn. Res.},
  FJOURNAL = {Journal of Machine Learning Research (JMLR)},
    VOLUME = {18},
      YEAR = {2017},
     PAGES = {Paper No. 165, 29},
      ISSN = {1532-4435,1533-7928},
   MRCLASS = {68T05 (49N90 62M45)},
  MRNUMBER = {3813814},
}

@inproceedings{
lu2019understanding,
title={Understanding and Improving Transformer From a Multi-Particle Dynamic System Point of View},
author={Yiping Lu and others},
booktitle={ICLR 2020 Workshop on Integration of Deep Neural Models and Differential Equations},
year={2019}
}

@article {burger2025analysis,
    AUTHOR = {Burger, Martin and Kabri, Samira and Korolev, Yury and Roith,
              Tim and Weigand, Lukas},
     TITLE = {Analysis of mean-field models arising from self-attention
              dynamics in transformer architectures with layer
              normalization},
   JOURNAL = {Philos. Trans. Roy. Soc. A},
  FJOURNAL = {Philosophical Transactions of the Royal Society A.
              Mathematical, Physical and Engineering Sciences},
    VOLUME = {383},
      YEAR = {2025},
    NUMBER = {2298},
     PAGES = {Paper No. 20240233, 48},
      ISSN = {1364-503X,1471-2962},
   MRCLASS = {58D20 (49Q22 65M75)},
  MRNUMBER = {4935785},
}

@article {geshkovski2025mathematical,
    AUTHOR = {Geshkovski, Borjan and Letrouit, Cyril and Polyanskiy, Yury
              and Rigollet, Philippe},
     TITLE = {A mathematical perspective on transformers},
   JOURNAL = {Bull. Amer. Math. Soc. (N.S.)},
  FJOURNAL = {American Mathematical Society. Bulletin. New Series},
    VOLUME = {62},
      YEAR = {2025},
    NUMBER = {3},
     PAGES = {427--479},
      ISSN = {0273-0979,1088-9485},
   MRCLASS = {34D05 (34C15 34D06 35Q83 52C17 60K35 82C22)},
  MRNUMBER = {4926874}
}

@inproceedings{sanderSinkformers2022,
  title = {Sinkformers: {{Transformers}} with {{Doubly Stochastic Attention}}},
  booktitle = {Proceedings of {{The}} 25th {{International Conference}} on {{Artificial Intelligence}} and {{Statistics}}},
  author = {Sander, M. E. and Ablin, P. and Blondel, M. and Peyr{\'e}, G.},
  year = {2022},
  pages = {3515--3530},
}

@article{geshkovski2023emergence,
  title={The emergence of clusters in self-attention dynamics},
  author={Geshkovski, Borjan and Letrouit, Cyril and Polyanskiy, Yury and Rigollet, Philippe},
  journal={Advances in Neural Information Processing Systems},
  volume={36},
  pages={57026--57037},
  year={2023}
}

@book {brezis2011functional,
    AUTHOR = {Brezis, Haim},
     TITLE = {Functional analysis, {S}obolev spaces and partial differential
              equations},
    SERIES = {Universitext},
 PUBLISHER = {Springer, New York},
      YEAR = {2011},
     PAGES = {xiv+599},
      ISBN = {978-0-387-70913-0},
   MRCLASS = {35-01 (46-01 46E35 46N20 47F05)},
  MRNUMBER = {2759829},
MRREVIEWER = {Vicen\c tiu\ D.\ R\u adulescu},
}

@article {alvarez2025constructive,
    AUTHOR = {\'Alvarez-L\'opez, Antonio and Geshkovski, Borjan and
              Ruiz-Balet, Dom\`enec},
     TITLE = {Constructive approximate transport maps with normalizing
              flows},
   JOURNAL = {Appl. Math. Optim.},
  FJOURNAL = {Applied Mathematics and Optimization},
    VOLUME = {92},
      YEAR = {2025},
    NUMBER = {2},
     PAGES = {Paper No. 33},
      ISSN = {0095-4616,1432-0606},
   MRCLASS = {35Q49 (46N10 68T07 93C15 93C20)},
  MRNUMBER = {4954519}
}

@article {wiener1932tauberian,
    AUTHOR = {Wiener, Norbert},
     TITLE = {Tauberian theorems},
   JOURNAL = {Ann. of Math. (2)},
  FJOURNAL = {Annals of Mathematics. Second Series},
    VOLUME = {33},
      YEAR = {1932},
    NUMBER = {1},
     PAGES = {1--100},
      ISSN = {0003-486X,1939-8980},
   MRCLASS = {99-04},
  MRNUMBER = {1503035}
}

@article {raissi2019physics,
    AUTHOR = {Raissi, M. and Perdikaris, P. and Karniadakis, G. E.},
     TITLE = {Physics-informed neural networks: a deep learning framework
              for solving forward and inverse problems involving nonlinear
              partial differential equations},
   JOURNAL = {J. Comput. Phys.},
  FJOURNAL = {Journal of Computational Physics},
    VOLUME = {378},
      YEAR = {2019},
     PAGES = {686--707},
      ISSN = {0021-9991,1090-2716},
   MRCLASS = {65M70 (68T05)},
  MRNUMBER = {3881695}
}

@article {cybenko1989approximation,
    AUTHOR = {Cybenko, G.},
     TITLE = {Approximation by superpositions of a sigmoidal function},
   JOURNAL = {Math. Control Signals Systems},
  FJOURNAL = {Mathematics of Control, Signals, and Systems},
    VOLUME = {2},
      YEAR = {1989},
    NUMBER = {4},
     PAGES = {303--314},
      ISSN = {0932-4194,1435-568X},
   MRCLASS = {41A30 (92A09 92A25)},
  MRNUMBER = {1015670},
MRREVIEWER = {A.\ Haimovici}
}

@book {wiener1948cybernetics,
    AUTHOR = {Wiener, Norbert},
     TITLE = {Cybernetics, or {C}ontrol and {C}ommunication in the {A}nimal
              and the {M}achine},
    SERIES = {Actualit\'es Scientifiques et Industrielles [Current
              Scientific and Industrial Topics]},
    VOLUME = {No. 1053},
 PUBLISHER = {Hermann \& Cie, Paris; The Technology Press, Cambridge, MA;
              John Wiley \& Sons, Inc., New York},
      YEAR = {1948},
     PAGES = {194},
   MRCLASS = {60.0X},
  MRNUMBER = {25096},
MRREVIEWER = {J.\ L.\ Doob},
}

@article{alvlop2023optimized, 
author = {\'Alvarez-L\'opez, Antonio and Orive-Illera, Rafael and Zuazua, Enrique}, 
title = {Cluster-Based Classification with Neural ODEs via Control}, journal = {Journal of Machine Learning}, 
year = {2025}, 
volume = {4}, 
number = {2}, pages = {128--156}, 
issn = {2790-2048}
}

@article {DomenecNormalisingFlows,
    AUTHOR = {Ruiz-Balet, Dom\`enec and Zuazua, Enrique},
     TITLE = {Control of neural transport for normalising flows},
   JOURNAL = {J. Math. Pures Appl. (9)},
  FJOURNAL = {Journal de Math\'ematiques Pures et Appliqu\'ees. Neuvi\`eme
              S\'erie},
    VOLUME = {181},
      YEAR = {2024},
     PAGES = {58--90},
      ISSN = {0021-7824,1776-3371},
   MRCLASS = {93B05 (35F10 49Q22 68T07)},
  MRNUMBER = {4668814},
MRREVIEWER = {Borjan\ Geshkovski}
}

@article {geshkovski2022turnpike,
    AUTHOR = {Geshkovski, Borjan and Zuazua, Enrique},
     TITLE = {Turnpike in optimal control of {PDE}s, {R}es{N}ets, and
              beyond},
   JOURNAL = {Acta Numer.},
  FJOURNAL = {Acta Numerica},
    VOLUME = {31},
      YEAR = {2022},
     PAGES = {135--263},
      ISSN = {0962-4929,1474-0508},
   MRCLASS = {49K20 (35Q49 91B62)},
  MRNUMBER = {4436586},
MRREVIEWER = {Martin\ Gugat}
}

@article{alcalde2025clustering,
author = {Alcalde, Albert and Fantuzzi, Giovanni and Zuazua, Enrique},
title = {Clustering in Pure-Attention Hardmax Transformers and Its Role in Sentiment Analysis},
journal = {SIAM Journal on Mathematics of Data Science},
volume = {7},
number = {3},
pages = {1367-1393},
year = {2025}
}

@article{alcalde2025exact,
  title={Exact Sequence Classification with Hardmax Transformers},
  author={Alcalde, Albert and Fantuzzi, Giovanni and Zuazua, Enrique},
  journal={arXiv preprint arXiv:2502.02270},
  year={2025}
}

@article{zhou2022towards,
  title={Towards understanding the condensation of neural networks at initial training},
  author={Zhou, Hanxu and Qixuan, Zhou and Luo, Tao and Zhang, Yaoyu and Xu, Zhi-Qin},
  journal={Advances in Neural Information Processing Systems},
  volume={35},
  pages={2184--2196},
  year={2022}
}

@book{zuazuaober,
author={Zuazua, Enrique},
title={Oberwolfach Seminars: Machine Learning with PDE and Control Perspectives},
    year={2026, to appear},
      publisher = {Birkh\"auser}
      }

@article{lorenzothysez,
  title={HYCO: Hybrid-Cooperative Learning for Data-Driven PDE Modeling},
  author={Liverani, Lorenzo and Steynberg, Matthys and Zuazua, Enrique},
  journal={arXiv preprint arXiv:2509.14123},
  year={2025}
}

@article{wang2023approximatea,
  title = {Approximate and Weighted Data Reconstruction Attack in Federated Learning},
  author = {Wang, Ziqi and Song, Yongcun and Zuazua, Enrique},
  year = {2023},
  journal = {arXiv preprint 	arXiv:2308.06822}
}

@article{song2024fedadmm,
  author = {Yongcun Song and Ziqi Wang and Enrique Zuazua},
  title = {FedADMM-InSa: An Inexact and Self-Adaptive ADMM for Federated Learning},
  journal = {Neural Networks},
  volume = {181},
  pages = {106772},
  year = {2025}
}

@inproceedings{he2016deep,
  title = {Deep Residual Learning for Image Recognition},
  booktitle = {Proceedings of the {{IEEE}} Conference on Computer Vision and Pattern Recognition},
  author = {He, Kaiming and Zhang, Xiangyu and Ren, Shaoqing and Sun, Jian},
  year = {2016},
  pages = {770--778}
}

@article{alvlop2024interplay,
  author = {Antonio \'Alvarez-L\'opez and Arselane Hadj Slimane and Enrique Zuazua},
  title = {Interplay between depth and width for interpolation in neural ODEs},
  journal = {Neural Networks},
  pages = {106640},
  year = {2024}
}

@article {bach2017,
    AUTHOR = {Bach, Francis},
     TITLE = {Breaking the curse of dimensionality with convex neural
              networks},
   JOURNAL = {J. Mach. Learn. Res.},
  FJOURNAL = {Journal of Machine Learning Research (JMLR)},
    VOLUME = {18},
      YEAR = {2017},
     PAGES = {Paper No. 19, 53},
      ISSN = {1532-4435,1533-7928},
   MRCLASS = {92B20 (62G05 90C25)},
  MRNUMBER = {3634886},
}

@article {devore_hanin_petrova_2021,
    AUTHOR = {DeVore, Ronald and Hanin, Boris and Petrova, Guergana},
     TITLE = {Neural network approximation},
   JOURNAL = {Acta Numer.},
  FJOURNAL = {Acta Numerica},
    VOLUME = {30},
      YEAR = {2021},
     PAGES = {327--444},
      ISSN = {0962-4929,1474-0508},
   MRCLASS = {68T07 (05C90)},
  MRNUMBER = {4298220},
MRREVIEWER = {Claudio\ Carvalhaes}
}

@article {domenecNODES,
    AUTHOR = {Ruiz-Balet, Dom\`enec and Zuazua, Enrique},
     TITLE = {Neural {ODE} control for classification, approximation, and
              transport},
   JOURNAL = {SIAM Rev.},
  FJOURNAL = {SIAM Review},
    VOLUME = {65},
      YEAR = {2023},
    NUMBER = {3},
     PAGES = {735--773},
      ISSN = {0036-1445,1095-7200},
   MRCLASS = {34H05 (35Q49 37N35 49N90 68T07 93C10)},
  MRNUMBER = {4624336},
MRREVIEWER = {Morteza\ Pakdaman}
}

@article {fernandez2000cost,
    AUTHOR = {Fern\'andez-Cara, Enrique and Zuazua, Enrique},
     TITLE = {The cost of approximate controllability for heat equations:
              the linear case},
   JOURNAL = {Adv. Differential Equations},
  FJOURNAL = {Advances in Differential Equations},
    VOLUME = {5},
      YEAR = {2000},
    NUMBER = {4-6},
     PAGES = {465--514},
      ISSN = {1079-9389},
   MRCLASS = {93B05 (35B37 35K05 93C20)},
  MRNUMBER = {1750109},
MRREVIEWER = {Werner\ Horn},
}

@article{FCZ,
  author = {Enrique Fern\'andez-Cara and Enrique Zuazua},
  title = {Control theory: History, mathematical achievements and perspectives},
  journal = {Bol. Soc. Esp. Mat. Apl.},
  volume = {26},
  pages = {79-140},
  year = {2003}
}

@article{Hernandez2024Classification,
  title={Deep neural networks: Multi-classification and universal approximation},
  author={Hern{\'a}ndez, Mart{\'\i}n and Zuazua, Enrique},
  journal={arXiv preprint arXiv:2409.06555},
  year={2024}
}

@article{liu2025moments,
  title={Moments, Time-Inversion and Source Identification for the Heat Equation},
  author={Liu, Kang and Zuazua, Enrique},
  journal={arXiv preprint arXiv:2507.02677},
  year={2025}
}

@article {liyau,
    AUTHOR = {Li, Peter and Yau, Shing-Tung},
     TITLE = {On the parabolic kernel of the {S}chr\"odinger operator},
   JOURNAL = {Acta Math.},
  FJOURNAL = {Acta Mathematica},
    VOLUME = {156},
      YEAR = {1986},
    NUMBER = {3-4},
     PAGES = {153--201},
      ISSN = {0001-5962,1871-2509},
   MRCLASS = {58G11 (35J10)},
  MRNUMBER = {834612},
MRREVIEWER = {Harold\ Donnelly}
}

@article{li2024universal,
  title={Universal approximation of dynamical systems by semi-autonomous neural odes and applications},
  author={Li, Ziqian and Liu, Kang and Liverani, Lorenzo and Zuazua, Enrique},
  journal={arXiv preprint arXiv:2407.17092},
  year={2024}
}

@article{liu2024fedgame,
  title={A Potential Game Perspective in Federated Learning},
  author={Liu, Kang and Wang, Ziqi and Zuazua, Enrique},
  journal={arXiv preprint arXiv:2411.11793},
  year={2024}
}

@article {liu2024representation,
    AUTHOR = {Liu, Kang and Zuazua, Enrique},
     TITLE = {Representation and regression problems in neural networks:
              relaxation, generalization, and numerics},
   JOURNAL = {Math. Models Methods Appl. Sci.},
  FJOURNAL = {Mathematical Models and Methods in Applied Sciences},
    VOLUME = {35},
      YEAR = {2025},
    NUMBER = {6},
     PAGES = {1471--1521},
      ISSN = {0218-2025,1793-6314},
   MRCLASS = {65K10 (68T07 68T09 90C06 90C26)},
  MRNUMBER = {4902962}
}

@article {papamakarios_normalizing_2021,
    AUTHOR = {Papamakarios, George and Nalisnick, Eric and Rezende, Danilo
              Jimenez and Mohamed, Shakir and Lakshminarayanan, Balaji},
     TITLE = {Normalizing flows for probabilistic modeling and inference},
   JOURNAL = {J. Mach. Learn. Res.},
  FJOURNAL = {Journal of Machine Learning Research (JMLR)},
    VOLUME = {22},
      YEAR = {2021},
     PAGES = {Paper No. 57, 64},
      ISSN = {1532-4435,1533-7928},
   MRCLASS = {62G99 (28A33 60E10)},
  MRNUMBER = {4253750},
}

@article {fisher1975spline,
    AUTHOR = {Fisher, S. D. and Jerome, J. W.},
     TITLE = {Spline solutions to {$L\sp{1}$} extremal problems in one and
              several variables},
   JOURNAL = {J. Approximation Theory},
  FJOURNAL = {Journal of Approximation Theory},
    VOLUME = {13},
      YEAR = {1975},
     PAGES = {73--83},
      ISSN = {0021-9045,1096-0430},
   MRCLASS = {41A65},
  MRNUMBER = {361577},
MRREVIEWER = {W.\ W.\ Breckner}
}

@article{zuazuanews,
  author = {Enrique Zuazua},
  title = {Control and Machine Learning},
  journal = {SIAM News},
  volume = {55},
  number = {8},
  month = {October},
  year = {2022}
}
\end{document}